\documentclass[reqno,12pt]{amsart}

\newtheorem{Theorem}{Theorem}[section]
\newtheorem{Proposition}[Theorem]{Proposition}

\newtheorem{Remark}[Theorem]{Remark}
\numberwithin{equation}{section}

\allowdisplaybreaks

\author[Michael J.\ Schlosser]{Michael J.\ Schlosser$^*$}
\address{Fakult\"at f\"ur Mathematik, Universit\"at Wien,
Nordbergstra{\ss}e 15, A-1090 Vienna, Austria}
\email{michael.schlosser@univie.ac.at}
\urladdr{http://www.mat.univie.ac.at/{\textasciitilde}schlosse}
\thanks{$^*$Partly supported by FWF Austrian Science Fund
grants \hbox{P17563-N13}, and S9607 (the second is part
of the Austrian National Research Network
``Analytic Combinatorics and Probabilistic Number Theory'').}
\date{August 29, 2006}


\title[Multilateral inversion of $A_r$, $C_r$ and $D_r$ series]
{Multilateral inversion of $\boldsymbol A_{\boldsymbol r}$,
$\boldsymbol C_{\boldsymbol r}$ and $\boldsymbol D_{\boldsymbol r}$\\
basic hypergeometric series}

\begin{document}

\begin{abstract}
In [Electron.\ J.\ Combin.\ \textbf{10} (2003), \#R10], the author
presented a new basic hypergeometric matrix inverse with applications
to bilateral basic hypergeometric series. This matrix inversion result
was directly extracted from an instance of Bailey's very-well-poised
${}_6\psi_6$ summation theorem, and involves two infinite matrices
which are not lower-triangular. The present paper features three
different multivariable generalizations of the above result.
These are extracted from Gustafson's $A_r$ and $C_r$ extensions
and of the author's recent $A_r$ extension of Bailey's
${}_6\psi_6$ summation formula. By combining
these new multidimensional matrix inverses with $A_r$ and $D_r$
extensions of Jackson's ${}_8\phi_7$ summation theorem three
balanced very-well-poised ${}_8\psi_8$ summation theorems associated
with the root systems $A_r$ and $C_r$ are derived.
\end{abstract}

\maketitle

\section{Introduction}\label{sec0}

In \cite[Th.~3.1]{S2}, the author presented the following matrix inverse:

Let $|q|<1$, and $a$, $b$ and $c$ be indeterminates. The infinite matrices
$(f_{nk})_{n,k\in\mathbb Z}$ and $(g_{kl})_{k,l\in\mathbb Z}$ are {\em
inverses} of each other where
\begin{subequations}\label{bmi}
\begin{align}\label{bmf}\notag
f_{nk}=\frac{(aq/b,bq/a,aq/c,cq/a,bq,q/b,cq,q/c)_{\infty}}
{(q,q,aq,q/a,aq/bc,bcq/a,cq/b,bq/c)_{\infty}}&\\\times
\frac{(1-bcq^{2n}/a)}{(1-bc/a)}\frac{(b)_{n+k}\,(a/c)_{k-n}}
{(cq)_{n+k}\,(aq/b)_{k-n}}&
\end{align}
and
\begin{equation}\label{bmg}
g_{kl}=\frac{(1-aq^{2k})}{(1-a)}\frac{(c)_{k+l}\,(a/b)_{k-l}}
{(bq)_{k+l}\,(aq/c)_{k-l}}\,q^{k-l}.
\end{equation}
\end{subequations}
(The notation is explained in Section~\ref{secpre}.)

This result was directly extracted from an instance of
Bailey's~\cite[Eq.~(4.7)]{B} very-well-poised ${}_6\psi_6$
summation formula,
\begin{multline}\label{66gl}
{}_6\psi_6\!\left[\begin{matrix}q\sqrt{a},-q\sqrt{a},b,c,d,e\\
\sqrt{a},-\sqrt{a},aq/b,aq/c,aq/d,aq/e\end{matrix}\,;q,
\frac{a^2q}{bcde}\right]\\
=\frac {(q,aq,q/a,aq/bc,aq/bd,aq/be,aq/cd,aq/ce,aq/de)_{\infty}}
{(aq/b,aq/c,aq/d,aq/e,q/b,q/c,q/d,q/e,a^2q/bcde)_{\infty}},
\end{multline}
where $|a^2q/bcde|<1$ (cf.\ \cite[Eq.~(5.3.1)]{GR}).

If we let $c\to a$ in \eqref{bmi}, we obtain a matrix inverse
found by Bressoud~\cite{Br} which he directly extracted from the
terminating very-well-poised ${}_6\phi_5$ summation
(a special case of \eqref{66gl}). If, after letting $c\to a$,
we additionally let $a\to 0$, we obtain Andrews'~\cite[Lemma~3]{A1}
``Bailey transform matrices'', a matrix inversion underlying
the powerful Bailey lemma. While Bressoud's matrix inverse
underlies Andrews' WP-Bailey lemma \cite{A2} (WP stands for
``well-poised'') which generalizes the classical Bailey lemma,
the ``bilateral'' matrix inverse \eqref{bmi} underlies the
BWP-Bailey lemma, a bilateral generalization of the WP-Bailey
lemma, see \cite{S5}. 

In \cite{S2}, several applications of \eqref{bmi} to bilateral basic
hypergeometric series were given. One of them included a new
very-well-poised ${}_8\psi_8$ summation formula,
see Proposition~\ref{88n} in this paper.

Here we apply part of the analysis of \cite{S2} to multiple sums.
In fact, by appropriately specializing Gustafson's $A_r$ and $C_r$
${}_6\psi_6$ summations~\cite{G1,G2}, and an $A_r$ ${}_6\psi_6$ summation
by the author~\cite{S4}, we derive three multidimensional extensions
of the bilateral matrix inverse \eqref{bmi} and deduce three
multidimensional ${}_8\psi_8$ summations, associated with the
root systems of type $A_r$ and $C_r$, as applications.
These are obtained via {\em multidimensional inverse relations},
applied to $A_r$ and $D_r$ extensions of Jackson's terminating balanced
very-well-poised ${}_8\phi_7$ summations, taken from \cite{M1,S1,S4}.

Our paper is organized as follows. In Section~\ref{secpre},
we first cover some preliminaries on basic hypergeometric series.
In the same section, we also explain some facts we need
on multidimensional basic hypergeometric series associated with
root systems. We list several multi-sum identities explicitly there
for easy reference.
Section~\ref{secmmi} is devoted to multidimensional matrix inversions.
In particular, we give three new explicit multilateral matrix inverses,
which are directly extracted from corresponding multivariate
${}_6\psi_6$ summation formulae. Our applications, see Section~\ref{sec88},
include three balanced very-well-poised ${}_8\psi_8$ summation formulae,
two of them associated with the root system $A_r$,
the third with the root system $C_r$. These new multivariate
${}_8\psi_8$ summations comprise, via specialization and analytic
continuation, corresponding multivariate ${}_8\phi_7$ and
${}_6\psi_6$ summation formulae. Finally, we show in the Appendix how
an incorrect application of multidimensional inverse relations
leads to a false result, namely a divergent $D_r$ very-well-poised
${}_6\psi_6$ summation (which however remains true for $r=1$,
or whenever the series terminates).

\section{Preliminaries}\label{secpre}

\subsection{Basic hypergeometric series}

We recall some standard notation for basic
hypergeometric series (cf.~\cite{GR}), and then turn to some
selected identities.

Let $q$ be a complex number such that $0<|q|<1$. We define the
{\em $q$-shifted factorial}\/ for all integers $k$ by
\begin{equation*}
(a)_{\infty}:=\prod_{j=0}^{\infty}(1-aq^j)\qquad\text{and}\qquad
(a)_k:=\frac{(a)_{\infty}}{(aq^k)_{\infty}}.
\end{equation*}
For brevity, we employ the condensed notation
\begin{equation*}
(a_1,\ldots,a_m)_k:= (a_1)_k\dots(a_m)_k
\end{equation*}
where $k$ is an integer or infinity. Further, we utilize
\begin{equation}\label{defhyp}
{}_s\phi_{s-1}\!\left[\begin{matrix}a_1,a_2,\dots,a_s\\
b_1,b_2,\dots,b_{s-1}\end{matrix}\,;q,z\right]:=
\sum _{k=0} ^{\infty}\frac {(a_1,a_2,\dots,a_s)_k}
{(q,b_1,\dots,b_{s-1})_k}z^k,
\end{equation}
and
\begin{equation}\label{defhypb}
{}_s\psi_s\!\left[\begin{matrix}a_1,a_2,\dots,a_s\\
b_1,b_2,\dots,b_s\end{matrix}\,;q,z\right]:=
\sum _{k=-\infty} ^{\infty}\frac {(a_1,a_2,\dots,a_s)_k}
{(b_1,b_2,\dots,b_s)_k}z^k,
\end{equation}
to denote the {\em basic hypergeometric ${}_s\phi_{s-1}$ series},
and the {\em bilateral basic hypergeometric ${}_s\psi_s$ series},
respectively. In \eqref{defhyp} or \eqref{defhypb}, $a_1,\dots,a_s$ are
called the {\em upper parameters}, $b_1,\dots,b_s$ the
{\em lower parameters}, $z$ is the {\em argument}, and 
$q$ the {\em base} of the series.
See \cite[p.~5 and p.~137]{GR} for the criteria
of when these series terminate, or, if not, when they converge. 
 
The classical theory of basic hypergeometric series contains
numerous summation and transformation formulae
involving ${}_s\phi_{s-1}$ or ${}_s\psi_s$ series.
Many of these summation theorems require
that the parameters satisfy the condition of being
either balanced and/or very-well-poised.
An ${}_s\phi_{s-1}$ basic hypergeometric series is called
{\em balanced} if $b_1\cdots b_{s-1}=a_1\cdots a_sq$ and $z=q$.
An ${}_s\phi_{s-1}$ series is {\em well-poised} if
$a_1q=a_2b_1=\cdots=a_sb_{s-1}$. An ${}_s\phi_{s-1}$ basic
hypergeometric series is called {\em very-well-poised}
if it is well-poised and if $a_2=-a_3=q\sqrt{a_1}$.
Note that the factor
\begin{equation}\label{vwp}
\frac {1-a_1q^{2k}}{1-a_1}
\end{equation}
appears in a very-well-poised series.
The parameter $a_1$ is usually referred to as the
{\em special parameter} of such a series.
Similarly, a bilateral ${}_s\psi_s$ basic hypergeometric series is
well-poised if $a_1b_1=a_2b_2\cdots=a_sb_s$ and very-well-poised if,
in addition, $a_1=-a_2=qb_1=-qb_2$. Further, we call a bilateral
${}_s\psi_s$ basic hypergeometric series balanced if
$b_1\cdots b_s=a_1\cdots a_sq^2$ and $z=q$.

A standard reference for basic hypergeometric series
is Gasper and Rahman's texts~\cite{GR}.
In our computations in Sections~\ref{secmmi} and \ref{sec88}
we frequently use some elementary identities of
$q$-shifted factorials, listed in \cite[Appendix~I]{GR}.

One of the most important theorems in the theory of basic
hypergeometric series is 
F.~H.~Jackson's~\cite{J} terminating balanced very-well-poised
${}_8\phi_7$ summation (cf.\ \cite[Eq.~(2.6.2)]{GR}):
\begin{multline}\label{87gl}
{}_8\phi_7\!\left[\begin{matrix}a,\,q\sqrt{a},-q\sqrt{a},b,c,d,
a^2q^{1+n}/bcd,q^{-n}\\
\sqrt{a},-\sqrt{a},aq/b,aq/c,aq/d,bcdq^{-n}/a,aq^{1+n}\end{matrix}\,;q,
q\right]\\
=\frac {(aq,aq/bc,aq/bd,aq/cd)_n}
{(aq/b,aq/c,aq/d,aq/bcd)_n}.
\end{multline}
A combinatorial proof of the {\em elliptic} extension of \eqref{87gl},
namely of Frenkel and Turaev's~\cite{FT}
${}_{10}V_9$ summation, which degenerates
to a combinatorial proof of \eqref{87gl} in the {\em trigonometric}
special case, has recently been given in \cite{S3}.

In \cite[Thm.~4.1]{S2}, Jackson's summation \eqref{87gl} was utilized,
in conjunction with the bilateral matrix inverse \eqref{bmi}, to derive 
the following balanced very-well-poised ${}_8\psi_8$ summation formula:

\begin{Proposition}\label{88n}
Let $a$, $b$, $c$ and $d$ be indeterminates,
let $k$ be an arbitrary integer and $M$ a nonnegative integer. Then
\begin{multline}\label{88ngl}
{}_8\psi_8\!\left[\begin{matrix}q\sqrt{a},-q\sqrt{a},b,c,dq^k,
aq^{-k}/c,aq^{1+M}/b,aq^{-M}/d\\
\sqrt{a},-\sqrt{a},aq/b,aq/c,aq^{1-k}/d,
cq^{1+k},bq^{-M},dq^{1+M}\end{matrix}\,;q,q\right]\\
=\frac {(aq/bc,cq/b,dq,dq/a)_M}
{(cdq/a,dq/c,q/b,aq/b)_M}
\frac{(cd/a,bd/a,cq,cq/a,dq^{1+M}/b,q^{-M})_k}
{(q,cq/b,d/a,d,bcq^{-M}/a,cdq^{1+M}/a)_k}\\\times
\frac{(q,q,aq,q/a,cdq/a,aq/cd,cq/d,dq/c)_{\infty}}
{(cq,q/c,dq,q/d,cq/a,aq/c,dq/a,aq/d)_{\infty}}.
\end{multline}
\end{Proposition}
Note that two of the upper parameters of the ${}_8\psi_8$ series
in \eqref{88ngl} (namely $b$ and $aq^{1+M}/b$) differ
multiplicatively from corresponding lower parameters by $q^M$,
(namely $bq^{-M}$ and $aq/b$, respectively) a nonnegative integral
power of $q$.

One can also derive (or verify) \eqref{88ngl} by adequately specializing
M.~Jackson's~\cite[Eq.~(2.2)]{Ja} transformation
formula for a very-well-poised $_8\psi_8$ series into a sum of two
(multiples of) $_8\phi_7$ series (cf.\ \cite[Eq.~(5.6.2)]{GR}):
\begin{multline}\label{88tgl}
{}_8\psi_8\!\left[\begin{matrix}q\sqrt{a},-q\sqrt{a},b,c,d,e,f,g\\
\sqrt{a},-\sqrt{a},aq/b,aq/c,aq/d,aq/e,aq/f,aq/g\end{matrix}\,;q,
\frac{a^3q^2}{bcdefg}\right]\\
=\frac{(q,aq,q/a,d,d/a,bq/c,bq/e,bq/f,bq/g,aq/bc,aq/be,aq/bf,
aq/bg)_{\infty}}
{(q/b,q/c,q/e,q/f,q/g,aq/b,aq/c,aq/e,aq/f,aq/g,d/b,bd/a,
b^2q/a)_{\infty}}\\\times
{}_8\phi_7\!\left[\begin{matrix}b^2/a,\,qb/\sqrt{a},-qb/\sqrt{a},
bc/a,bd/a,be/a,bf/a,bg/a\\
b/\sqrt{a},-b/\sqrt{a},bq/c,bq/d,bq/e,bq/f,bq/g\end{matrix}\,;q,
\frac{a^3q^2}{bcdefg}\right]\\
+\frac{(q,aq,q/a,b,b/a,dq/c,dq/e,dq/f,dq/g,aq/cd,aq/de,aq/df,
aq/dg)_{\infty}}
{(q/c,q/d,q/e,q/f,q/g,aq/c,aq/d,aq/e,aq/f,aq/g,b/d,bd/a,
d^2q/a)_{\infty}}\\\times
{}_8\phi_7\!\left[\begin{matrix}d^2/a,\,qd/\sqrt{a},-qd/\sqrt{a},
bd/a,cd/a,de/a,df/a,dg/a\\
d/\sqrt{a},-d/\sqrt{a},dq/b,dq/c,dq/e,dq/f,dq/g\end{matrix}\,;q,
\frac{a^3q^2}{bcdefg}\right],
\end{multline}
where $|a^3q^2/bcdefg|<1$, for convergence.
In particular, substituting $d\mapsto dq^k$, and then letting
$e\to aq^{-k}/b$, $f\to aq^{1+M}/b$, $g\to aq^{-M}/d$,
the coefficient of the first $_8\phi_7$ series on the right-hand side
becomes zero (as it contains $(q^{-M})_\infty$),
while the second $_8\phi_7$ series can be summed
by an application of Jackson's terminating $_8\phi_7$
summation in \eqref{87gl}.
This way of establishing \eqref{88ngl} works so far only in
the classical one-dimensional case, as no multiple series
extensions of \eqref{88tgl} are yet known.
Our multivariate extensions of Proposition~\ref{88n}
in Section~\ref{sec88} of this paper, see
Theorems~\ref{a88s} and \ref{c88s} (obtained by suitable
extensions of the analysis applied in \cite{S2}), which we find
attractive by themselves, can be understood as a first step in
the quest of finding multivariate extensions of the very-well-poised
$_8\psi_8$ transformation formula \eqref{88tgl}, or of even more general
transformations.

Two special cases of Proposition~\ref{88n} are worth pointing out:
\begin{enumerate}
\item
If $c\to a$ (or $c\to q^{-k}$), then the bilateral series in \eqref{88ngl}
gets truncated from below and from above so that the sum is finite. By
a polynomial argument, $q^M$ can then be replaced by any complex number.
If we replace it by $bc/aq$, perform the substitutions
$d\mapsto dq^{-k}$, and finally replace $k$ by $n$,
we obtain exactly Jackson's terminating balanced
very-well-poised ${}_8\phi_7$ summation in \eqref{87gl}.

\item
If, in \eqref{88ngl}, we let $M\to\infty$, perform the substitution
$d\mapsto dq^{-k}$, and rewrite the products of the form $(x)_k$ as
$(x)_\infty/(xq^k)_\infty$, we can apply analytic continuation
to replace $q^k$ by $a/ce$ (in order to relax the integrality condition
of $k$) where $e$ is a new complex parameter. We then obtain exactly 
Bailey's very-well-poised ${}_6\psi_6$ summation in \eqref{66gl}.
\end{enumerate}

\subsection{Multidimensional basic hypergeometric series associated
with root systems}\label{secmult}

$A_r$ (or, equivalently, $U(r+1)$) hypergeometric series were motivated
by the work of Biedenharn, Holman, and Louck~\cite{HBL}
in theoretical physics.
The theory of $A_r$ basic hypergeometric series
(or ``multiple basic hypergeometric series associated with the
root system $A_r$'', or ``associated with the unitary group $U(r+1)$''),
analogous to the classical theory of one-dimensional series,
has been developed originally by R.\ A.\ Gustafson, S.\ C.\ Milne,
and their co-workers, and later others (see \cite{G1,G2,M1,M2,ML,Ro}
for a very small selection of papers in this area).
Notably, several higher-dimensional extensions have been derived
(in each case) for the $q$-binomial theorem, $q$-Chu--Vandermonde
summation, $q$-Pfaff--Saalsch\"utz summation, Jackson's ${}_8\phi_7$
summation, Bailey's ${}_{10}\phi_9$ transformation, and
other important summation and transformation theorems.
See \cite{M3} for a survey on some of the main results and
techniques from the theory of $A_r$ basic hypergeometric series.
Multiple basic hypergeometric series associated with other roots systems
than $A_r$ have been first defined by Gustafson~\cite{G2}
who succeeded in giving several multivariable extensions
of Bailey's ${}_6\psi_6$ summation.
In particular, some important results for $C_r$ and $D_r$ basic
hypergeometric series have been derived in \cite{Bh,BS,DG,G2,LM,ML,S1}
(-- again, this is a very imcomplete listing).

We note the conventions for naming our series as 
$A_r$, $C_r$ or $D_r$ basic hypergeometric series.
We consider multiple series of the form 
$\sum_{k_1,k_2,\dots,k_r} S({\mathbf k})$,
where ${\mathbf k}=(k_1,\dots,k_r)$,
which reduce to classical basic hypergeometric series when $r=1$.
We call such a multiple basic hypergeometric 
series {\em balanced} if it reduces to 
a balanced series when $r=1$. Well-poised and very-well-poised
series are defined similarly. 

Further, such a multiple series 
is called a $C_r$ basic hypergeometric series if the summand
$S({\mathbf k})$ contains the factor
\begin{equation}\label{cnvandy}
\prod_{1\le i<j\le r}\frac{x_iq^{k_i}-x_jq^{k_j}}{x_i-x_j}
\prod_{1\le i\le j\le r}\frac{1-x_ix_jq^{k_i+k_j}}{1-x_ix_j}
\end{equation} 
Note that when $r=1$, \eqref{cnvandy} reduces to 
\begin{equation*}
\frac {1-x_1^2q^{2k_1}}{1-x_1^2},
\end{equation*}
which is \eqref{vwp} with $x_1^2$ acting like the special parameter
of a very-well poised series. In our statements of $C_r$ theorems,
we set $x_i\mapsto {\sqrt{a}}x_i$ for $i=1,\dots,r$, and make
similar changes to other parameters in $S({\mathbf k})$.
This is done in order to follow the classical notation in \cite{GR} as
closely as possible. A typical example of a $C_r$ basic hypergeometric
series is the left-hand side of \eqref{cr87gl}.

$D_r$ multiple basic series are closely related to $C_r$ series.
Instead of \eqref{cnvandy}, $S({\mathbf k})$ only has the following factors:
\begin{equation}\label{dnvandy}
\prod_{1\le i<j\le r}\frac{(x_iq^{k_i}-x_jq^{k_j})
(1-x_ix_jq^{k_i+k_j})}
{(x_i-x_j)(1-x_ix_j)}.
\end{equation}
A typical example is the left-hand side of \eqref{dr87gl}
(with $x_i\mapsto {\sqrt{cd}}x_i$ for $i=1,\dots,r$).

$A_r$ basic hypergeometric series only have 
\begin{equation}\label{anvandy}
\prod_{1\le i<j\le r}\frac{x_iq^{k_i}-x_jq^{k_j}}{x_i-x_j}
\end{equation} 
as a factor of $S({\mathbf k})$.
Typical examples are the left-hand sides of \eqref{ar87gl}
and \eqref{arv87gl}.
A reason for naming these
series as $A_r$, $C_r$ or $D_r$ series is that \eqref{anvandy},
\eqref{dnvandy}, and \eqref{cnvandy} are closely associated with 
the product side of the Weyl denominator formulae
for the respective root systems, see \cite{Bh,RS,St}.

For compact notation, we usually write
\begin{equation*}
|{\mathbf k}|:=k_1+\dots+k_r,\quad\text{where}\quad
{\mathbf k}=(k_1,\dots,k_r),
\end{equation*}
and
\begin{equation*}
C:=c_1\cdots c_r,\quad E:=e_1\cdots e_r,\quad\text{etc.}
\end{equation*}

We now list several multivariable extensions of Jackson's ${}_8\phi_7$
summation \eqref{87gl}. The first identity is taken from
\cite[Thm.~6.17]{M1}.

\begin{Proposition}[({\sc Milne}) An $A_r$ terminating balanced
very-well-poised ${}_8\phi_7$ summation formula]\label{ar87}
Let $a$, $b$, $c_1,\dots,c_r$, $d$ and $x_1,\dots,x_r$ be
indeterminate, and let $M$ be a nonnegative integer. Then
\begin{multline}\label{ar87gl}
\underset{0\le|\mathbf k|\le M}{\sum_{k_1,\dots,k_r\ge 0}}
\prod_{1\le i<j\le r}\frac {x_iq^{k_i}-x_jq^{k_j}}{x_i-x_j}
\prod_{i=1}^r\frac{1-ax_iq^{k_i+|{\mathbf k}|}}{1-ax_i}
\prod_{i,j=1}^r\frac{(c_jx_i/x_j)_{k_i}}{(qx_i/x_j)_{k_i}}\\\times
\prod_{i=1}^r\frac{(ax_i)_{|{\mathbf k}|}\,
(dx_i,a^2x_iq^{1+M}/bCd)_{k_i}}
{(ax_iq/c_i)_{|{\mathbf k}|}\,(ax_iq/b,ax_iq^{1+M})_{k_i}}\cdot
\frac{(b,q^{-M})_{|{\mathbf k}|}}{(aq/d,bCdq^{-M}/a)_{|{\mathbf k}|}}\,
q^{|{\mathbf k}|}\\
=\frac{(aq/bd,aq/Cd)_M}{(aq/d,aq/bCd)_M}\,
\prod_{i=1}^r\frac{(ax_iq,ax_iq/bc_i)_M}
{(ax_iq/b,ax_iq/c_i)_M}.
\end{multline}
\end{Proposition}

The following identity was recently obtained in \cite[Eq.~(4.3)]{S4}.

\begin{Proposition}[({\sc S.}) An $A_r$ terminating balanced 
very-well-poised ${}_8\phi_7$ summation formula]\label{arv87}
Let $a$, $b$, $c_1,\dots,c_r$, $d$ and $x_1,\dots,x_r$ be
indeterminate, and let $M$ be a nonnegative integer. Then
\begin{multline}\label{arv87gl}
\underset{0\le|{\mathbf k}|\le M}{\sum_{k_1,\dots,k_r\ge 0}}
\frac{(1-aq^{2|{\mathbf k}|})}{(1-a)}
\prod_{1\le i<j\le r}\frac{x_iq^{k_i}-x_jq^{k_j}}{x_i-x_j}
\prod_{i,j=1}^r\frac{(c_jx_i/x_j)_{k_i}}{(qx_i/x_j)_{k_i}}\\\times
\prod_{i=1}^r\frac{(aq/Cx_id)_{|{\mathbf k}|-k_i}\,
(b/x_i)_{|{\mathbf k}|}\,(dx_i)_{k_i}}
{(b/x_i)_{|{\mathbf k}|-k_i}\,(ac_iq/Cx_id)_{|{\mathbf k}|}\,
(ax_iq/b)_{k_i}}\cdot
\frac{(a,a^2q^{1+M}/bCd,q^{-M})_{|{\mathbf k}|}}
{(aq/C,bCdq^{-M}/a,aq^{1+M})_{|{\mathbf k}|}}\,
q^{|{\mathbf k}|}\\=
\frac{(aq,aq/bd)_M}{(aq/C,aq/bCd)_M}
\prod_{i=1}^r\frac{(aq/Cx_id,ax_iq/bc_i)_M}{(ax_iq/b,ac_iq/Cx_id)_M}.
\end{multline}
\end{Proposition}

The following identity was derived in \cite[Thm.~4.1]{DG},
and, independently, in \cite[Thm.~6.13]{ML}.

\begin{Proposition}[({\sc Denis--Gustafson; Milne--Lilly})
A $C_r$ terminating balanced very-well-poised ${}_8\phi_7$
summation formula]\label{cr87}
Let $a$, $b$, $c$, $d$ and $x_1,\dots,x_r$ be indeterminate
and let $m_1,\dots,m_r$ be nonnegative integers. Then
\begin{multline}\label{cr87gl}
\underset{i=1,\dots,r}{\sum_{0\le k_i\le m_i}}
\prod_{1\le i<j\le r}\frac {x_iq^{k_i}-x_jq^{k_j}}{x_i-x_j}
\prod_{1\le i\le j\le r}\frac{1-ax_ix_jq^{k_i+k_j}}{1-ax_ix_j}
\prod_{i,j=1}^r\frac{(q^{-m_j}x_i/x_j,ax_ix_j)_{k_i}}
{(ax_ix_jq^{1+m_j},qx_i/x_j)_{k_i}}\\\times
\prod_{i=1}^r\frac{(bx_i,cx_i,dx_i,a^2x_iq^{1+|\mathbf m|}/bcd)_{k_i}}
{(ax_iq/b,ax_iq/c,ax_iq/d,bcdx_iq^{-|\mathbf m|}/a)_{k_i}}
\cdot q^{|{\mathbf k}|}\\
=\prod_{1\le i<j\le r}(ax_ix_jq)_{m_i+m_j}^{-1}
\prod_{i,j=1}^r(ax_ix_jq)_{m_i}\\\times
\frac{(aq/bc,aq/bd,aq/cd)_{|\mathbf m|}}
{\prod_{i=1}^r(ax_iq/b,ax_iq/c,ax_iq/d,aq^{1+|\mathbf m|-m_i}/bcdx_i)_{m_i}}.
\end{multline}
\end{Proposition}

The last extension of \eqref{87gl} we need
was derived in \cite[Thm.~5.17]{S1}.

\begin{Proposition}[({\sc S.}) A $D_r$ terminating balanced
very-well-poised ${}_8\phi_7$ summation formula]\label{dr87}
Let $a$, $b$, $c$, $c_1,\dots,c_r$, $d$, and $x_1,\dots,x_r$ be
indeterminate and let $M$ be a nonnegative integer. Then
\begin{multline}\label{dr87gl}
\underset{0\le|\mathbf k|\le M}{\sum_{k_1,\dots,k_r\ge 0}}
\prod_{1\le i<j\le r}\frac {x_iq^{k_i}-x_jq^{k_j}}{x_i-x_j}
\prod_{i=1}^r\frac{1-ax_iq^{k_i+|{\mathbf k}|}}{1-ax_i}
\prod_{i=1}^r\frac{(ax_i)_{|{\mathbf k}|}\,(aq/cdx_i)_{|\mathbf k|-k_i}}
{(ax_iq/c_i,ac_iq/cdx_i)_{|{\mathbf k}|}}\\\times
\prod_{1\le i<j\le r}(cdx_ix_j)_{k_i+k_j}^{-1}
\prod_{i,j=1}^r\frac{(c_jx_i/x_j,cdx_ix_j/c_j)_{k_i}}
{(qx_i/x_j)_{k_i}}\\\times
\frac{(b,a^2q^{1+M}/bcd,q^{-M})_{|\mathbf k|}}
{\prod_{i=1}^r(ax_iq/b,bcdx_iq^{-M}/a,ax_iq^{1+M})_{k_i}}\,
q^{|{\mathbf k}|}\\
=\prod_{i=1}^r\frac{(ax_iq,ax_iq/bc_i,ac_iq/bcdx_i,aq/cdx_i)_M}
{(aq/bcdx_i,ac_iq/cdx_i,ax_iq/c_i,ax_iq/b)_M}.
\end{multline}
\end{Proposition}

A closely related $D_r$ terminating balanced
very-well-poised ${}_8\phi_7$ summation, equivalent to
Proposition~\ref{dr87} by reversing summations
of the ``rectangular form'' of Proposition~\ref{dr87},
given in \cite[Thm.~5.6]{S1}, was derived by
Bhatnagar, see \cite[Thm.~7]{Bh}.

Next, we list several multivariable extensions of Bailey's ${}_6\psi_6$
summation in \eqref{66gl}. The first of these
was derived in \cite[Thm.~1.15]{G1}.

\begin{Proposition}[({\sc Gustafson}) An $A_r$ very-well-poised
${}_6\psi_6$ summation formula]\label{ar66}
Let $a$, $b$, $c_1,\dots,c_r$, $d$, $e_1,\dots,e_r$ and
$x_1,\dots,x_r$ be indeterminate. Then
\begin{multline}\label{ar66gl}
\sum_{k_1,\dots,k_r=-\infty}^{\infty}
\prod_{1\le i<j\le r}\frac {x_iq^{k_i}-x_jq^{k_j}}{x_i-x_j}
\prod_{i=1}^r\frac{1-ax_iq^{k_i+|{\mathbf k}|}}{1-ax_i}
\prod_{i,j=1}^r\frac{(c_jx_i/x_j)_{k_i}}
{(ax_iq/e_jx_j)_{k_i}}\\\times
\prod_{i=1}^r\frac{(e_ix_i)_{|{\mathbf k}|}\,(dx_i)_{k_i}}
{(ax_iq/c_i)_{|{\mathbf k}|}\,(ax_iq/b)_{k_i}}\cdot
\frac{(b)_{|{\mathbf k}|}}{(aq/d)_{|{\mathbf k}|}}
\left(\frac{a^{r+1}q}{bCdE}\right)^{|{\mathbf k}|}\\
=\frac{(aq/bd,a^rq/bE,aq/Cd)_{\infty}}
{(a^{r+1}q/bCdE,aq/d,q/b)_{\infty}}
\prod_{i,j=1}^r\frac{(ax_iq/c_ie_jx_j,qx_i/x_j)_{\infty}}
{(qx_i/c_ix_j,ax_iq/e_jx_j)_{\infty}}\\\times
\prod_{i=1}^r\frac{(ax_iq/bc_i,aq/de_ix_i,ax_iq,q/ax_i)_{\infty}}
{(ax_iq/b,ax_iq/c_i,q/dx_i,q/e_ix_i)_{\infty}},
\end{multline}
provided $|a^{r+1}q/bCdE|<1$.
\end{Proposition}

The following identity was recently obtained in \cite[Thm.~6.1]{S4}.

\begin{Proposition}[({\sc S.}) An $A_r$ very-well-poised
${}_6\psi_6$ summation formula]\label{arv66}
Let $a$, $b$, $c_1,\dots,c_r$, $d$, $e_1,\dots,e_r$ and
$x_1,\dots,x_r$ be indeterminate. Then
\begin{multline}\label{arv66gl}
\sum_{k_1,\dots,k_r=-\infty}^{\infty}
\frac{(1-aq^{2|{\mathbf k}|})}{(1-a)}
\prod_{1\le i<j\le r}\frac {x_iq^{k_i}-x_jq^{k_j}}{x_i-x_j}
\prod_{i,j=1}^r\frac{(c_jx_i/x_j)_{k_i}}{(ax_iq/e_jx_j)_{k_i}}
\\\times
\prod_{i=1}^r\frac{(aq/Cdx_i)_{|{\mathbf k}|-k_i}\
(bE/a^{r-1}e_ix_i)_{|{\mathbf k}|}\,(dx_i)_{k_i}}
{(bE/a^rx_i)_{|{\mathbf k}|-k_i}\,(ac_iq/Cdx_i)_{|{\mathbf k}|}\,
(ax_iq/b)_{k_i}}\cdot
\frac{(E/a^{r-1})_{|{\mathbf k}|}}{(aq/C)_{|{\mathbf k}|}}
\left(\frac{a^{r+1}q}{bCdE}\right)^{|{\mathbf k}|}\\=
\frac{(aq,q/a,aq/bd)_\infty}{(aq/C,a^{r+1}q/bCdE,a^{r-1}q/E)_\infty}
\prod_{i,j=1}^r\frac{(qx_i/x_j,ax_iq/c_ie_jx_j)_\infty}
{(qx_i/c_ix_j,ax_iq/e_jx_j)_\infty}\\\times
\prod_{i=1}^r\frac{(a^rx_iq/bE,aq/e_idx_i,aq/Cdx_i,ax_iq/bc_i)_\infty}
{(a^{r-1}e_ix_iq/bE,q/dx_i,ax_iq/b,ac_iq/Cdx_i)_\infty},
\end{multline}
provided $|aq^{r+1}/bCdE|<1$.
\end{Proposition}

The third extension of \eqref{66gl} we need is taken from
\cite[Thm.~5.1]{G2}.

\begin{Proposition}[({\sc Gustafson}) A $C_r$ very-well-poised
${}_6\psi_6$ summation formula]\label{cr66}
Let $a$, $b$, $c_1,\dots,c_r$, $d$, $e_1,\dots,e_r$ and
$x_1,\dots,x_r$ be indeterminate. Then
\begin{multline}\label{cr66gl}
\sum_{k_1,\dots,k_r=-\infty}^{\infty}
\prod_{1\le i<j\le r}\frac {x_iq^{k_i}-x_jq^{k_j}}{x_i-x_j}
\prod_{1\le i\le j\le r}\frac {1-ax_ix_jq^{k_i+k_j}}{1-ax_ix_j}\\\times
\prod_{i,j=1}^r\frac{(c_jx_i/x_j,e_jx_ix_j)_{k_i}}
{(ax_ix_jq/c_j,ax_iq/e_jx_j)_{k_i}}
\prod_{i=1}^r\frac{(bx_i,dx_i)_{k_i}}{(ax_iq/b,ax_iq/d)_{k_i}}\cdot
\left(\frac{a^{r+1}q}{bCdE}\right)^{|{\mathbf k}|}\\
=\prod_{1\le i<j\le r}(ax_ix_jq/c_ic_j,aq/e_ie_jx_ix_j)_{\infty}
\prod_{1\le i\le j\le r}(ax_ix_jq,q/ax_ix_j)_{\infty}\\\times
\frac{(aq/bd)_{\infty}}{(a^{r+1}q/bCdE)_{\infty}}
\prod_{i,j=1}^r\frac{(ax_iq/c_ie_jx_j,qx_i/x_j)_{\infty}}
{(ax_iq/e_jx_j,q/e_jx_ix_j,ax_ix_jq/c_i,qx_i/c_ix_j)_{\infty}}\\\times
\prod_{i=1}^r\frac{(ax_iq/bc_i,aq/be_ix_i,ax_iq/c_id,aq/de_ix_i)_{\infty}}
{(ax_iq/b,q/bx_i,ax_iq/d,q/dx_i)_{\infty}},
\end{multline}
provided $|a^{r+1}q/bCdE|<1$.
\end{Proposition}

Having listed several of the most fundamental summation formulae
of the theory of multidimensional basic hypergeometric series
associated with root systems,
we are now ready to turn to the derivation of new results.

\section{Multidimensional matrix inversions}\label{secmmi}

Let $\mathbb Z$ denote the set of integers.
In the following, we consider infinite $r$-dimensional matrices
$F=(f_{\mathbf n\mathbf k})_{\mathbf n,\mathbf k\in\mathbb Z^r}$ and
$G=(g_{\mathbf n\mathbf k})_{\mathbf n,\mathbf k\in\mathbb Z^r}$, and
infinite sequences $(a_{\mathbf n})_{\mathbf n\in\mathbb Z^r}$ and
$(b_{\mathbf n})_{\mathbf n\in\mathbb Z^r}$.

Clearly, $F$ is the {\em left-inverse} of $G$,
if and only if the following orthogonality
relation holds:
\begin{equation}\label{orthrel}
\sum_{\mathbf k\in\mathbb Z^r}f_{\mathbf n\mathbf k}g_{\mathbf k\mathbf l}
=\delta_{\mathbf n\mathbf l}\qquad\qquad
\text {for all}\quad \mathbf n,\mathbf l\in\mathbb Z^r.
\end{equation}
Further, $F$ is the {\em right-inverse} of $G$,
if and only if the following orthogonality
relation holds:
\begin{equation}\label{orthreld}
\sum_{\mathbf k\in\mathbb Z^r}g_{\mathbf n\mathbf k}f_{\mathbf k\mathbf l}
=\delta_{\mathbf n\mathbf l}\qquad\qquad
\text {for all}\quad \mathbf n,\mathbf l\in\mathbb Z^r.
\end{equation}

If $F$ is the left-inverse {\em and} the right-inverse of $G$ we simply
say that $F$ and $G$ are {\em mutually inverse} or
{\em inverses of each other}.

Note that in \eqref{orthrel} and \eqref{orthreld} we are {\em not}
requiring that the infinite $r$-dimensional matrices are lower-triangular
(which would mean that $f_{\mathbf n\mathbf k}=g_{\mathbf n\mathbf k}=0$
unless $\mathbf n\ge\mathbf k$, where by the latter we mean
$n_i\ge k_i$ for all $i=1,\dots,r$).
If they were lower-triangular, the multiple series on the left-hand
sides of \eqref{orthrel} and \eqref{orthreld} would be in fact finite
sums (and both relations must then hold at the same time).
In the general case, the sums are infinite.
If the summands of the infinite series involve complex numbers,
we require suitable convergence conditions to hold
(such as absolute convergence; for interchanging
double sums we further need uniform convergence).
Note that convergence of one of the sums does not
necessarily imply convergence of the other.

Now consider the following two equations (a.k.a.\ ``inverse relations''):
\begin{subequations}\label{invrel}
\begin{equation}\label{invrel1}
\sum_{\mathbf k\in\mathbb Z^r}f_{\mathbf n\mathbf k}a_{\mathbf k}
=b_{\mathbf n}\qquad\qquad\text{for all $\mathbf n$,}
\end{equation}
and
\begin{equation}\label{invrel2}
\sum_{\mathbf l\in\mathbb Z^r}g_{\mathbf k\mathbf l}b_{\mathbf l}
=a_{\mathbf k}\qquad\qquad\text{for all $\mathbf k$.}
\end{equation}
\end{subequations}
It is immediate from the orthogonality relations \eqref{orthrel}
and \eqref{orthreld} that if $F$ is the left-inverse of $G$,
the relation \eqref{invrel2} implies \eqref{invrel1},
while if $F$ is the right-inverse of $G$,
the relation \eqref{invrel1} implies \eqref{invrel2},
subject to convergence.

Similarly, one may consider another pair of equations, where one sums
over the {\em first} (instead of the second) multi-index of the matrix:
\begin{subequations}\label{rotinv}
\begin{equation}\label{rotinv1}
\sum_{\mathbf n\in\mathbb Z^r}f_{\mathbf n\mathbf k}a_{\mathbf n}
=b_{\mathbf k}\qquad\qquad\text{for all $\mathbf k$,}
\end{equation}
and
\begin{equation}\label{rotinv2}
\sum_{\mathbf k\in\mathbb Z^r}g_{\mathbf k\mathbf l}b_{\mathbf k}
=a_{\mathbf l}\qquad\qquad\text{for all $\mathbf l$.}
\end{equation}
\end{subequations}
Again, it is immediate from the orthogonality relations \eqref{orthrel}
and \eqref{orthreld} that if $F$ is the left-inverse of $G$,
the relation \eqref{rotinv1} implies \eqref{rotinv2},
while if $F$ is the right-inverse of $G$,
the relation \eqref{rotinv2} implies \eqref{rotinv1},
again subject to convergence. 

We are ready to state and prove three multidimensional matrix inverses,
all of them as consequences of corresponding multivariate ${}_6\psi_6$
summations which have been stated in Section~\ref{secpre}.

\begin{Theorem}[An $A_r$ multilateral matrix inverse]\label{armi}
Let $a$, $b$, $c_1,\dots,c_r$, and $x_1,\dots,x_r$ be
indeterminate. Then
$(f_{\mathbf n\mathbf k})_{\mathbf n,\mathbf k\in\mathbb Z^r}$ and
$(g_{\mathbf k\mathbf l})_{\mathbf k,\mathbf l\in\mathbb Z^r}$ are
{\em inverses} of each other where
\begin{subequations}
\begin{multline}\label{farmi}
f_{\mathbf n\mathbf k}=
\frac{(bq,q/b)_{\infty}}{(bq/C,Cq/b)_{\infty}}
\prod_{i,j=1}^r\frac{(qc_jx_i/x_j,qx_i/c_ix_j)_{\infty}}
{(qc_jx_i/c_ix_j,qx_i/x_j)_{\infty}}\\\times
\prod_{i=1}^r\frac{(ax_iq/b,bq/ax_i,ax_iq/c_i,c_iq/ax_i)_{\infty}}
{(ax_iq/bc_i,bc_iq/ax_i,ax_iq,q/ax_i)_{\infty}}\\\times
\prod_{1\le i<j\le r}\frac {c_iq^{n_i}/x_i-c_jq^{n_j}/x_j}{c_i/x_i-c_j/x_j}
\prod_{i=1}^r\frac{1-bc_iq^{n_i+|\mathbf n|}/ax_i}{1-bc_i/ax_i}\\\times
(b)_{|\mathbf n|+|\mathbf k|}
\prod_{i,j=1}^r\frac 1{(qc_jx_i/x_j)_{n_j+k_i}}
\prod_{i=1}^r\frac{(ax_i/c_i)_{|\mathbf k|-n_i}}
{(ax_iq/b)_{k_i-|\mathbf n|}}
\end{multline}
and
\begin{multline}\label{garmi}
g_{\mathbf k \mathbf l}=
\prod_{1\le i<j\le r}\frac {x_iq^{k_i}-x_jq^{k_j}}{x_i-x_j}
\prod_{i=1}^r\frac{1-ax_iq^{k_i+|\mathbf k|}}{1-ax_i}\\\times
\frac 1{(bq)_{|\mathbf k|+|\mathbf l|}}\,
\prod_{i,j=1}^r(c_jx_i/x_j)_{k_i+l_j}
\prod_{i=1}^r\frac{(ax_i/b)_{k_i-|\mathbf l|}}
{(ax_iq/c_i)_{|\mathbf k|-l_i}}
\cdot q^{|\mathbf k|-r|\mathbf l|}.
\end{multline}
\end{subequations}
\end{Theorem}

\begin{proof}
We show that the inverse matrices \eqref{farmi}/\eqref{garmi}
satisfy the orthogonality relation \eqref{orthrel}.
(An analogous computation reveals that the inverse matrices
\eqref{farmi}/\eqref{garmi} also satisfy the dual orthogonality
relation \eqref{orthreld}.) Writing out the sum
$\sum_{\mathbf k\in\mathbb Z^r}f_{\mathbf n\mathbf k}g_{\mathbf k\mathbf l}$
with the above choices of $f_{\mathbf n\mathbf k}$ and
$g_{\mathbf k\mathbf l}$ we observe that the multiple series can be
summed by an application of the $A_r$ very-well-poised ${}_6\psi_6$
summation in Proposition~\ref{ar66}. The specializations needed there are
\begin{equation}\label{subs}
b\mapsto bq^{|\mathbf n|},\quad\; c_i\mapsto c_iq^{l_i},\quad\;
d\mapsto aq^{-|\mathbf l|}/b, \quad\; e_i\mapsto aq^{-n_i}/c_i,
\quad\; i=1,\dots,r.
\end{equation}
The summation formula gives us a product containing the factors
\begin{equation}\label{ardnl}
(q^{1+|\mathbf l|-|\mathbf n|})_\infty
\prod_{i,j=1}^r(q^{1+n_j-l_i}c_jx_i/c_ix_j)_\infty.
\end{equation}
Since \eqref{ardnl} vanishes for all $r$-tuples of integers $\mathbf n$
and $\mathbf l$ with $\mathbf n\neq\mathbf l$, we can simplify the
product (setting $\mathbf n=\mathbf l$, the only non-zero case) and
readily determine that the sum indeed boils down to
$\delta_{\mathbf n\mathbf l}$. The details are as follows:
\begin{multline}\label{auxid}
\sum_{\mathbf k\in\mathbb Z^r}f_{\mathbf n\mathbf k}g_{\mathbf k\mathbf l}=
\sum_{k_1,\dots,k_r=-\infty}^\infty
\frac{(bq,q/b)_{\infty}}{(bq/C,Cq/b)_{\infty}}
\prod_{i,j=1}^r\frac{(qc_jx_i/x_j,qx_i/c_ix_j)_{\infty}}
{(qc_jx_i/c_ix_j,qx_i/x_j)_{\infty}}\\\times
\prod_{i=1}^r\frac{(ax_iq/b,bq/ax_i,ax_iq/c_i,c_iq/ax_i)_{\infty}}
{(ax_iq/bc_i,bc_iq/ax_i,ax_iq,q/ax_i)_{\infty}}\\\times
\prod_{1\le i<j\le r}\frac {c_iq^{n_i}/x_i-c_jq^{n_j}/x_j}{c_i/x_i-c_j/x_j}
\prod_{i=1}^r\frac{1-bc_iq^{n_i+|\mathbf n|}/ax_i}{1-bc_i/ax_i}\\\times
(b)_{|\mathbf n|+|\mathbf k|}
\prod_{i,j=1}^r\frac 1{(qc_jx_i/x_j)_{n_j+k_i}}
\prod_{i=1}^r\frac{(ax_i/c_i)_{|\mathbf k|-n_i}}
{(ax_iq/b)_{k_i-|\mathbf n|}}\\\times
\prod_{1\le i<j\le r}\frac {x_iq^{k_i}-x_jq^{k_j}}{x_i-x_j}
\prod_{i=1}^r\frac{1-ax_iq^{k_i+|\mathbf k|}}{1-ax_i}\\\times
\frac 1{(bq)_{|\mathbf k|+|\mathbf l|}}\,
\prod_{i,j=1}^r(c_jx_i/x_j)_{k_i+l_j}
\prod_{i=1}^r\frac{(ax_i/b)_{k_i-|\mathbf l|}}
{(ax_iq/c_i)_{|\mathbf k|-l_i}}
\cdot q^{|\mathbf k|-r|\mathbf l|}\\=
\frac{(bq,q/b)_{\infty}}{(bq/C,Cq/b)_{\infty}}
\prod_{i,j=1}^r\frac{(qc_jx_i/x_j,qx_i/c_ix_j)_{\infty}}
{(qc_jx_i/c_ix_j,qx_i/x_j)_{\infty}}\\\times
\prod_{i=1}^r\frac{(ax_iq/b,bq/ax_i,ax_iq/c_i,c_iq/ax_i)_{\infty}}
{(ax_iq/bc_i,bc_iq/ax_i,ax_iq,q/ax_i)_{\infty}}\\\times
\prod_{1\le i<j\le r}\frac {c_iq^{n_i}/x_i-c_jq^{n_j}/x_j}{c_i/x_i-c_j/x_j}
\prod_{i=1}^r\frac{1-bc_iq^{n_i+|\mathbf n|}/ax_i}{1-bc_i/ax_i}\\\times
(b)_{|\mathbf n|}
\prod_{i,j=1}^r\frac 1{(qc_jx_i/x_j)_{n_j}}
\prod_{i=1}^r\frac{(ax_i/c_i)_{-n_i}}
{(ax_iq/b)_{-|\mathbf n|}}\\\times
\frac 1{(bq)_{|\mathbf l|}}\,
\prod_{i,j=1}^r(c_jx_i/x_j)_{l_j}
\prod_{i=1}^r\frac{(ax_i/b)_{-|\mathbf l|}}
{(ax_iq/c_i)_{-l_i}}
\cdot q^{-r|\mathbf l|}\\\times
\sum_{k_1,\dots,k_r=-\infty}^\infty
\prod_{1\le i<j\le r}\frac {x_iq^{k_i}-x_jq^{k_j}}{x_i-x_j}
\prod_{i=1}^r\frac{1-ax_iq^{k_i+|\mathbf k|}}{1-ax_i}
\prod_{i,j=1}^r\frac{(q^{l_j}c_jx_i/x_j)_{k_i}}
{(q^{1+n_j}c_jx_i/x_j)_{k_i}}\\\times
\prod_{i=1}^r\frac{(ax_iq^{-n_i}/c_i)_{|\mathbf k|}\,
(ax_iq^{-|\mathbf l|}/b)_{k_i}}
{(ax_iq^{1-l_i}/c_i)_{|\mathbf k|}\,(ax_iq^{1-|\mathbf n|}/b)_{k_i}}
\cdot\frac{(bq^{|\mathbf n|})_{|\mathbf k|}}
{(bq^{1+|\mathbf l|})_{|\mathbf k|}}\,q^{|\mathbf k|}\\
=\frac{(bq,q/b)_{\infty}}{(bq/C,Cq/b)_{\infty}}
\prod_{i,j=1}^r\frac{(qc_jx_i/x_j,qx_i/c_ix_j)_{\infty}}
{(qc_jx_i/c_ix_j,qx_i/x_j)_{\infty}}\\\times
\prod_{i=1}^r\frac{(ax_iq/b,bq/ax_i,ax_iq/c_i,c_iq/ax_i)_{\infty}}
{(ax_iq/bc_i,bc_iq/ax_i,ax_iq,q/ax_i)_{\infty}}\\\times
\prod_{1\le i<j\le r}\frac {c_iq^{n_i}/x_i-c_jq^{n_j}/x_j}{c_i/x_i-c_j/x_j}
\prod_{i=1}^r\frac{1-bc_iq^{n_i+|\mathbf n|}/ax_i}{1-bc_i/ax_i}\\\times
\prod_{i,j=1}^r\frac{(c_jx_i/x_j)_{l_j}}{(qc_jx_i/x_j)_{n_j}}
\prod_{i=1}^r\frac{(ax_i/c_i)_{-n_i}\,(ax_i/b)_{-|\mathbf l|}}
{(ax_iq/c_i)_{-l_i}\,(ax_iq/b)_{-|\mathbf n|}}\cdot
\frac{(b)_{|\mathbf n|}}{(bq)_{|\mathbf l|}}q^{-r|\mathbf l|}\\\times
\frac{(q^{1+|\mathbf l|-|\mathbf n|},Cq/b,bq/C)_{\infty}}
{(q,bq^{1+|\mathbf l|},q^{1-|\mathbf n|}/b)_{\infty}}
\prod_{i,j=1}^r\frac{(q^{1+n_j-l_i}c_jx_i/c_ix_j,qx_i/x_j)_{\infty}}
{(q^{1-l_i}x_i/c_ix_j,q^{1+n_j}c_jx_i/x_j)_{\infty}}\\\times
\prod_{i=1}^r\frac{(ax_iq^{1-l_i-|\mathbf n|}/bc_i,
bc_iq^{1+n_i+|\mathbf l|}/ax_i,ax_iq,q/ax_i)_{\infty}}
{(ax_iq^{1-|\mathbf n|}/b,ax_iq^{1-l_i}/c_i,bq^{1+|\mathbf l|}/ax_i,
c_iq^{1+n_i}/ax_i)_{\infty}}.
\end{multline}
Now we set $\mathbf n=\mathbf l$, apply several elementary
identities from \cite[App.~I]{GR} and apply the $n\mapsto r$,
$x_i\mapsto c_iq^{-n_i}/x_i$, $y_i\mapsto n_i$, $i=1,\dots,r$, case
of \cite[Lem.~3.12]{M2}, specifically
\begin{multline}\label{mid}
\prod_{i,j=1}^r(q^{1+n_j-n_i}c_jx_i/c_ix_j)_{n_i-n_j}=
(-1)^{(r-1)|\mathbf n|}\,q^{\binom{|\mathbf n|}2-
r\sum_{i=1}^r\binom{n_i}2}\\\times
\prod_{i=1}^r\left(\frac{c_i}{x_i}\right)^{|{\mathbf n}|-rn_i}
\prod_{1\le i<j\le r}\frac{c_i/x_i-c_j/x_j}{c_iq^{n_i}/x_i-c_jq^{n_j}/x_j},
\end{multline}
to transform the last expression obtained in \eqref{auxid} to
$\delta_{\mathbf n\mathbf l}$.
\end{proof}

\begin{Remark}\rm\label{rem1}
The $c_i\to 1$, $i=1,\dots,r$, case of Theorem~\ref{armi}
can be reduced to Milne's~\cite[Thm.~3.41]{M2} $A_r$ extension of
Bressoud's matrix inverse \cite{Br}. In particular, since
$1/(q)_{n+k}=0$ for $n+k<0$, and $(1)_{k+l}=0$ for $k+l>0$,
the orthogonality relation \eqref{orthrel} then reduces to
\begin{equation*}
\underset{i=1,\dots,r}{\sum_{-n_i\le k_i\le -l_i}}
f_{\mathbf n\mathbf k}g_{\mathbf k\mathbf l}=\delta_{\mathbf n\mathbf l},
\end{equation*}
i.e. (after replacing $\mathbf k$ by $-\mathbf k$),
\begin{equation*}
\underset{i=1,\dots,r}{\sum_{l_i\le k_i\le n_i}}
f_{\mathbf n,-\mathbf k}g_{-\mathbf k,\mathbf l}=\delta_{\mathbf n\mathbf l}.
\end{equation*}
The $r$-dimensional matrices
$(f_{\mathbf n,-\mathbf k})_{\mathbf n,\mathbf k\in\mathbb Z^r}$ and
$(g_{-\mathbf k,\mathbf l})_{\mathbf k,\mathbf l\in\mathbb Z^r}$
are thus mutually inverse {\em lower-triangular} matrices.
\end{Remark}

\begin{Theorem}[Another $A_r$ multilateral matrix inverse]\label{arvmi}
Let $a$, $b$, $c_1,\dots,c_r$ and $x_1,\dots,x_r$ be
indeterminate. Then
$(f_{\mathbf n\mathbf k})_{\mathbf n,\mathbf k\in\mathbb Z^r}$ and
$(g_{\mathbf k\mathbf l})_{\mathbf k,\mathbf l\in\mathbb Z^r}$ are
{\em inverses} of each other where
\begin{subequations}
\begin{multline}\label{farvmi}
f_{\mathbf n\mathbf k}=
\frac{(aq/C,Cq/a)_{\infty}}{(aq,q/a)_{\infty}}
\prod_{i,j=1}^r\frac{(qc_jx_i/x_j,qx_i/c_ix_j)_{\infty}}
{(qc_jx_i/c_ix_j,qx_i/x_j)_{\infty}}\\\times
\prod_{i=1}^r\frac{(ax_iq/b,bq/ax_i,bc_iq/Cx_i,Cx_iq/bc_i)_{\infty}}
{(ax_iq/bc_i,bc_iq/ax_i,Cx_iq/b,bq/Cx_i)_{\infty}}\\\times
\prod_{1\le i<j\le r}\frac {c_iq^{n_i}/x_i-c_jq^{n_j}/x_j}{c_i/x_i-c_j/x_j}
\prod_{i=1}^r\frac{1-bc_iq^{n_i+|\mathbf n|}/ax_i}{1-bc_i/ax_i}\\\times
(a/C)_{|\mathbf k|-|\mathbf n|}
\prod_{i,j=1}^r\frac 1{(qc_jx_i/x_j)_{n_j+k_i}}
\prod_{i=1}^r\frac{(bc_i/Cx_i)_{n_i+|\mathbf k|}}
{(ax_iq/b)_{k_i-|\mathbf n|}}
\end{multline}
and
\begin{multline}\label{garvmi}
g_{\mathbf k \mathbf l}=
\frac{(1-aq^{2|\mathbf k|})}{(1-a)}
\prod_{1\le i<j\le r}\frac {x_iq^{k_i}-x_jq^{k_j}}{x_i-x_j}
\prod_{i=1}^r\frac{1-bq^{|\mathbf k|-k_i}/Cx_i}{1-b/Cx_i}\\\times
\frac 1{(aq/C)_{|\mathbf k|-|\mathbf l|}}\,
\prod_{i,j=1}^r(c_jx_i/x_j)_{k_i+l_j}
\prod_{i=1}^r\frac{(ax_i/b)_{k_i-|\mathbf l|}}
{(bc_iq/Cx_i)_{|\mathbf k|+l_i}}
\cdot q^{|\mathbf k|-r|\mathbf l|}.
\end{multline}
\end{subequations}
\end{Theorem}

\begin{proof}
The proof is similar to the proof of Theorem~\ref{armi}
but utilizes Proposition~\ref{arv66} in addition to Proposition~\ref{ar66}.

We first show that the inverse matrices \eqref{farvmi}/\eqref{garvmi}
satisfy the orthogonality relation \eqref{orthrel}.
Writing out the sum
$\sum_{\mathbf k\in\mathbb Z^r}f_{\mathbf n\mathbf k}g_{\mathbf k\mathbf l}$
with the above choices of $f_{\mathbf n\mathbf k}$ and
$g_{\mathbf k\mathbf l}$ we observe that the multiple series can be
summed by an application of the $A_r$ very-well-poised ${}_6\psi_6$
summation in Proposition~\ref{arv66}. The specializations needed there
are exactly the same as in Equation~\eqref{subs}. Again, the summation
formula leads to a product containing the factors in \eqref{ardnl}.
We can thus simplify the product (setting $\mathbf n=\mathbf l$,
the only non-zero case) and readily determine, by applying several
elementary identities for $q$-shifted factorials, including
\eqref{mid}, that the sum indeed
boils down to $\delta_{\mathbf n\mathbf l}$.

An analogous computation reveals that the inverse matrices
\eqref{farvmi}/\eqref{garvmi} also satisfy the dual orthogonality relation
\eqref{orthreld}. Writing out the sum
$\sum_{\mathbf k\in\mathbb Z^r}g_{\mathbf n\mathbf k}f_{\mathbf k\mathbf l}$
with the above choices of $g_{\mathbf n\mathbf k}$ and
$f_{\mathbf k\mathbf l}$ we observe that the multiple series can be
summed by an application of the $A_r$ very-well-poised ${}_6\psi_6$
summation in Proposition~\ref{ar66}. The specializations needed there are
\begin{align*}
a&\mapsto b/a,& b&\mapsto Cq^{-|\mathbf n|}/a,&
d&\mapsto bq^{|\mathbf l|}/C, &\\
c_i&\mapsto c_iq^{n_i},& e_i&\mapsto bq^{-l_i}/ac_i,&
x_i&\mapsto c_i/x_i,&
i=1,\dots,r.
\end{align*}
The summation formula gives us a product containing the factors
\begin{equation}\label{ardnl1}
(q^{1-|\mathbf l|+|\mathbf n|})_\infty
\prod_{i,j=1}^r(q^{1+l_j-n_i}c_ix_j/c_jx_i)_\infty.
\end{equation}
Since \eqref{ardnl1} vanishes for all $r$-tuples of integers $\mathbf n$
and $\mathbf l$ with $\mathbf n\neq\mathbf l$, we can simplify the
product (setting $\mathbf n=\mathbf l$, the only non-zero case) and
readily determine, by applying several
elementary identities for $q$-shifted factorials, 
that the sum indeed boils down to
$\delta_{\mathbf n\mathbf l}$. We omit the details,
being similar to those as in the proof of Theorem~\ref{armi}.
\end{proof}

\begin{Remark}\rm
The $c_i\to 1$, $i=1,\dots,r$, case of Theorem~\ref{arvmi}
can be reduced to the author's~\cite[Cor.~3.2]{S4} $A_r$ extension of
Bressoud's matrix inverse \cite{Br}, which can be also obtained by
specializing Bhatnagar and Milne's matrix inverse \cite[Thm.~3.48]{BM}.
As in Remark~\ref{rem1}, the $r$-dimensional matrices
$(f_{\mathbf n,-\mathbf k})_{\mathbf n,\mathbf k\in\mathbb Z^r}$ and
$(g_{-\mathbf k,\mathbf l})_{\mathbf k,\mathbf l\in\mathbb Z^r}$
are then mutually inverse {\em lower-triangular} matrices.
\end{Remark}

The following matrix inverse serves as a bridge between $C_r$ series
and $D_r$ series.

\begin{Theorem}[A $C_r$/$D_r$ multilateral matrix inverse]\label{crmi}
Let $a$, $b$, $c_1,\dots,c_r$ and $x_1,\dots,x_r$ be
indeterminate. Then
$(f_{\mathbf n\mathbf k})_{\mathbf n,\mathbf k\in\mathbb Z^r}$
is the {\em left-inverse} of
$(g_{\mathbf k\mathbf l})_{\mathbf k,\mathbf l\in\mathbb Z^r}$
where
\begin{subequations}
\begin{multline}\label{fcrmi}
f_{\mathbf n\mathbf k}=
\prod_{i=1}^r\frac{(ax_iq/b,bq/ax_i,bx_iq,q/bx_i)_{\infty}}
{(ax_iq/bc_i,bc_iq/ax_i,bx_iq/c_i,c_iq/bx_i)_{\infty}}\\\times
\frac{\prod_{i,j=1}^r(qc_jx_i/x_j,qx_i/c_ix_j,
ax_ix_jq/c_i,c_jq/ax_ix_j)_{\infty}}
{\prod_{i,j=1}^r(qc_jx_i/c_ix_j,qx_i/x_j)_{\infty}\,
\prod_{1\le i\le j\le r}(ax_ix_jq,q/ax_ix_j)_{\infty}}\\\times
\prod\nolimits_{1\le i<j\le r}
(ax_ix_jq/c_ic_j,c_ic_jq/ax_ix_j)_{\infty}^{-1}\\\times
\prod_{1\le i<j\le r}\frac {c_iq^{n_i}/x_i-c_jq^{n_j}/x_j}{c_i/x_i-c_j/x_j}
\prod_{i,j=1}^r\frac{(ax_ix_j/c_j)_{k_i-n_j}}
{(qc_jx_i/x_j)_{k_i+n_j}}\\\times
\prod_{i=1}^r\frac{(1-bc_iq^{n_i+|\mathbf n|}/ax_i)
(1-bx_iq^{|\mathbf n|-n_i}/c_i)\,(bx_i)_{k_i+|\mathbf n|}}
{(1-bc_i/ax_i)(1-bx_i/c_i)\,(ax_iq/b)_{k_i-|\mathbf n|}}
\end{multline}
and
\begin{multline}\label{gcrmi}
g_{\mathbf k \mathbf l}=
\prod_{1\le i<j\le r}\frac {x_iq^{k_i}-x_jq^{k_j}}{x_i-x_j}
\prod_{1\le i\le j\le r}\frac {1-ax_ix_jq^{k_i+k_j}}{1-ax_ix_j}\\\times
\prod_{i,j=1}^r\frac{(c_jx_i/x_j)_{k_i+l_j}}
{(ax_ix_jq/c_j)_{k_i-l_j}}
\prod_{i=1}^r\frac{(ax_i/b)_{k_i-|\mathbf l|}}{(bx_iq)_{k_i+|\mathbf l|}}
\cdot q^{|\mathbf k|+(1-2r)|\mathbf l|}.
\end{multline}
\end{subequations}
\end{Theorem}

\begin{proof}
The proof is completely analogous to the proofs of
Theorems~\ref{armi} and \ref{arvmi}, with the only difference that
the orthogonality relation \eqref{orthrel} of the inverse matrices
\eqref{fcrmi}/\eqref{gcrmi} is now established by using
Proposition~\ref{cr66}. The specializations \eqref{subs} again
lead to a product containing the factors in \eqref{ardnl}.
\end{proof}

\begin{Remark}\label{rem}\rm
Theorem~\ref{crmi} states that
$(f_{\mathbf n\mathbf k})_{\mathbf n,\mathbf k\in\mathbb Z^r}$
is the {\em left-inverse} of
$(g_{\mathbf k\mathbf l})_{\mathbf k,\mathbf l\in\mathbb Z^r}$.
This is because the infinite multiple sum
$\sum_{\mathbf k\in\mathbb Z^r}f_{\mathbf n\mathbf k}g_{\mathbf k\mathbf l}$
converges for all $r=1,2\dots$, and evaluates to
$\delta_{\mathbf n\mathbf l}$.
On the contrary,
$(f_{\mathbf n\mathbf k})_{\mathbf n,\mathbf k\in\mathbb Z^r}$
is {\em not} the {\em right-inverse} of
$(g_{\mathbf k\mathbf l})_{\mathbf k,\mathbf l\in\mathbb Z^r}$
unless $r=1$. The infinite multiple sum
$\sum_{\mathbf k\in\mathbb Z^r}g_{\mathbf n\mathbf k}f_{\mathbf k\mathbf l}$
converges for all $r$ but for $r>1$ in general does {\em not}
evaluate to $\delta_{\mathbf n\mathbf l}$ (which we find quite surprising!).
In particular, we do {\em not} have a closed form evaluation for the
convergent multilateral sum
\begin{multline*}
\sum_{k_1,\dots,k_r=-\infty}^{\infty}
\prod_{1\le i<j\le r}\frac {c_iq^{k_i}/x_i-c_jq^{k_j}/x_j}{c_i/x_i-c_j/x_j}
\prod_{i=1}^r\frac{(1-bc_iq^{k_i+|\mathbf k|}/ax_i)
(1-bx_iq^{|\mathbf k|-k_i}/c_i)}
{(1-bc_i/ax_i)(1-bx_i/c_i)}\\\times
\prod_{i,j=1}^r\frac{(q^{n_i}c_jx_i/x_j,c_jq^{-n_i}/ax_ix_j)_{k_j}}
{(q^{1+l_i}c_jx_i/x_j,c_jq^{1-l_i}/ax_ix_j)_{k_j}}
\prod_{i=1}^r\frac{(bx_iq^{l_i},bq^{-l_i}/ax_i)_{|\mathbf k|}}
{(bx_iq^{1+n_i},bq^{1-n_i}/ax_i)_{|\mathbf k|}}\cdot q^{|\mathbf k|}.
\end{multline*}
To see what can go wrong, when one applies inverse relations
in an incorrect way, see the Appendix.
\end{Remark}

\begin{Remark}\rm
The $c_i\mapsto ax_i^2$, $i=1,\dots,r$, case of Theorem~\ref{crmi}
reduces to a multidimensional matrix inverse involving two
{\em lower-triangular} matrices being mutually inverse
(a $C_r/D_r$ extension of Bressoud's matrix inverse~\cite{Br}),
a result first derived in \cite[Thm.~5.10]{S1}.
The combination of this inversion with
Denis--Gustafson's/Milne--Lilly's~\cite{DG,ML} $C_r$ terminating
balanced very-well-poised ${}_8\phi_7$ summation, stated here as
Proposition~\ref{cr87}, led in \cite[Thm.~5.14]{S1} to a $D_r$
terminating balanced very-well-poised ${}_8\phi_7$ summation,
equivalent (by a polynomial argument)
to the $D_r$ summation stated here as Proposition~\ref{dr87}.
\end{Remark}

\section{Multivariable balanced very-well-poised
${}_8\psi_8$ summations}\label{sec88}

As applications of the multilateral matrix inverses of Section~\ref{secmmi},
we provide three multidimensional extensions of the ${}_8\psi_8$ summation
formula in Theorem~\ref{88n}.

\begin{Theorem}[An $A_r$ balanced very-well-poised
${}_8\psi_8$ summation formula]\label{a88s}
Let $a$, $b$, $c_1,\dots,c_r$, $d$ and $x_1,\dots,x_r$ be
indeterminate, let $k_1,\dots,k_r$ be integers, and let $M$ be a
nonnegative integer. Then
\begin{multline}\label{a88sgl}
\sum_{n_1,\dots,n_r=-\infty}^{\infty}
\prod_{1\le i<j\le r}\frac {x_iq^{n_i}-x_jq^{n_j}}{x_i-x_j}
\prod_{i=1}^r\frac{1-ax_iq^{n_i+|{\mathbf n}|}}{1-ax_i}
\prod_{i,j=1}^r\frac{(c_jx_i/x_j)_{n_i}}
{(q^{1+k_j}c_jx_i/x_j)_{n_i}}\\\times
\prod_{i=1}^r\frac{(ax_iq^{-k_i}/c_i)_{|{\mathbf n}|}\,
(bx_i,ax_iq^{-M}/d)_{n_i}}
{(ax_iq/c_i)_{|{\mathbf n}|}\,
(bx_iq^{-M},ax_iq^{1-|{\mathbf k}|}/d)_{n_i}}\cdot
\frac{(dq^{|\mathbf k|},aq^{1+M}/b)_{|{\mathbf n}|}}
{(dq^{1+M},aq/b)_{|{\mathbf n}|}}\,q^{|{\mathbf n}|}\\
=\prod_{i,j=1}^r\frac{(qc_jx_i/c_ix_j,qx_i/x_j)_{\infty}}
{(qc_jx_i/x_j,qx_i/c_ix_j)_{\infty}}
\prod_{i=1}^r\frac{(ax_iq,q/ax_i,ax_iq/c_id,c_idq/ax_i)_{\infty}}
{(ax_iq/c_i,c_iq/ax_i,ax_iq/d,dq/ax_i)_{\infty}}\\\times
\frac{(dq/C,Cq/d)_{\infty}}{(dq,q/d)_{\infty}}\,
\frac{(dq,aq/bC)_M}{(aq/b,dq/C)_M}\,
\prod_{i=1}^r\frac{(c_iq/bx_i,dq/ax_i)_M}
{(c_idq/ax_i,q/bx_i)_M}\prod_{i,j=1}^r\frac{(qc_jx_i/x_j)_{k_i}}
{(qc_jx_i/c_ix_j)_{k_i}}\\\times
\frac{(bd/a,q^{-M})_{|{\mathbf k}|}}
{(d,bCq^{-M}/a)_{|{\mathbf k}|}}
\prod_{i=1}^r\frac{(c_id/ax_i)_{|{\mathbf k}|}\,
(c_iq/ax_i,c_idq^{1+M}/bCx_i)_{k_i}}
{(d/ax_i)_{|{\mathbf k}|}\,(c_iq/bx_i,c_idq^{1+M}/ax_i)_{k_i}}.
\end{multline}
\end{Theorem}

\begin{proof}
We combine the multilateral matrix inverse in Theorem~\ref{armi}
with the $A_r$ extension of Jackson's terminating balanced
very-well-poised ${}_8\phi_7$ summation in Proposition~\ref{ar87},
using the inverse relations \eqref{rotinv}.
(Alternatively, we may also use the inverse relations \eqref{invrel},
with a similar analysis as in the proof of Theorem~\ref{av88s}.)

In particular, we have \eqref{rotinv2}, with
$(g_{\mathbf k\mathbf l})_{\mathbf k,\mathbf l\in\mathbb Z^r}$
as in \eqref{garmi},
\begin{multline*}
a_{\mathbf l}=\frac{(bq/d,bq/C)_M}{(bq,bq/Cd)_M}
\prod_{i=1}^r\frac{(ax_iq,ax_iq/c_id)_M}{(ax_iq/d,ax_iq/c_i)_M}
\prod_{i,j=1}^r(c_jx_i/x_j)_{l_j}\\\times
\frac{(bq^{1+M}/d)_{|\mathbf l|}}{(bq^{1+M},bq/d)_{|\mathbf l|}}
\prod_{i=1}^r\frac{(c_iq^{-M}/ax_i,c_id/ax_i)_{l_i}}
{(c_idq^{-M}/ax_i)_{l_i}\,(bq/ax_i)_{|\mathbf l|}}\\\times
(-1)^{(r-1)|\mathbf l|}a^{(1-r)|\mathbf l|}b^{r|\mathbf l|}
q^{r\binom{|\mathbf l|}2-\sum_{i=1}^r\binom{l_i}2}
\prod_{i=1}^rc_i^{-l_i}x_i^{l_i-|\mathbf l|},
\end{multline*}
and
\begin{equation*}
b_{\mathbf k}=\frac{(d,q^{-M})_{|\mathbf k|}}{(Cdq^{-M}/b)_{|\mathbf k|}}
\prod_{i=1}^r\frac{(abx_iq^{1+M}/Cd)_{k_i}\,(ax_i)_{|\mathbf k|}}
{(ax_iq/d,ax_iq^{1+M})_{k_i}}\prod_{i,j=1}^r\frac 1{(qx_i/x_j)_{k_i}},
\end{equation*}
by the $b\mapsto d$, $c_i\mapsto c_iq^{l_i}$,
$d\mapsto aq^{-|\mathbf l|}/b$, $i=1,\dots,r$, case of Proposition~\ref{ar87}.
Therefore we must have \eqref{rotinv1}, with
$(f_{\mathbf n\mathbf k})_{\mathbf n,\mathbf k\in\mathbb Z^r}$
as in \eqref{farmi}, and the above sequences $a_{\mathbf n}$
and $b_{\mathbf k}$. After simplifications and the substitutions
$a\mapsto d/a$, $b\mapsto d$, $d\mapsto bd/a$, $x_i\mapsto c_i/x_i$,
$i=1,\dots,r$, we arrive at \eqref{a88sgl}.
\end{proof}

\begin{proof}[Alternative proof of Theorem~\ref{a88s}]
We combine the multilateral matrix inverse in Theorem~\ref{arvmi}
with the $A_r$ extension of Jackson's terminating balanced
very-well-poised ${}_8\phi_7$ summation in Proposition~\ref{arv87},
using the inverse relations \eqref{rotinv}.

In particular, we have \eqref{rotinv2}, with
$(g_{\mathbf k\mathbf l})_{\mathbf k,\mathbf l\in\mathbb Z^r}$
as in \eqref{garvmi},
\begin{multline*}
a_{\mathbf l}=\frac{(aq,bq/d)_M}
{(aq/C,bq/Cd)_M}
\prod_{i=1}^r\frac{(bq/Cx_i,ax_iq/c_id)_M}
{(ax_iq/d,bc_iq/Cx_i)_M}
\prod_{i,j=1}^r(c_jx_i/x_j)_{l_j}\\\times
\frac{(bq^{1+M}/d,Cq^{-M}/a)_{|\mathbf l|}}{(bq/d)_{|\mathbf l|}}
\prod_{i=1}^r\frac{(c_id/ax_i)_{l_i}\,x_i^{|\mathbf l|}}
{(bc_iq^{1+M}/Cx_i,c_idq^{-M}/ax_i)_{l_i}\,(bx_iq/a)_{|\mathbf l|}}\\\times
(-1)^{(r-1)|\mathbf l|}
a^{(1-r)|\mathbf l|}b^{r|\mathbf l|}C^{-|\mathbf l|}
q^{(r-1)\binom{|\mathbf l|}2},
\end{multline*}
and
\begin{equation*}
b_{\mathbf k}=\frac{(a,abq^{1+M}/Cd,q^{-M})_{|\mathbf k|}}
{(Cdq^{-M}/b,aq^{1+M})_{|\mathbf k|}}
\prod_{i=1}^r\frac{(b/Cx_i)_{|\mathbf k|-k_i}\,(d/x_i)_{|\mathbf k|}}
{(d/x_i)_{|\mathbf k|-k_i}\,(ax_iq/d)_{k_i}}
\prod_{i,j=1}^r\frac 1{(qx_i/x_j)_{k_i}},
\end{equation*}
by the $b\mapsto d$, $c_i\mapsto c_iq^{l_i}$,
$d\mapsto aq^{-|\mathbf l|}/b$, $i=1,\dots,r$, case of Proposition~\ref{arv87}.
Therefore we must have \eqref{rotinv1}, with
$(f_{\mathbf n\mathbf k})_{\mathbf n,\mathbf k\in\mathbb Z^r}$
as in \eqref{farvmi}, and the above sequences $a_{\mathbf n}$
and $b_{\mathbf k}$. After simplifications and the substitutions
$a\mapsto d/a$, $b\mapsto d$, $d\mapsto bd/a$, $x_i\mapsto c_i/x_i$,
$i=1,\dots,r$, we arrive at \eqref{a88sgl}.
\end{proof}

\begin{Remark}\rm
Two special cases of Theorem~\ref{a88s} are of particular interest:
\begin{enumerate}
\item
If $c_i=q^{-k_i}$, for $i=1,\dots,r$, then the multilateral series in
\eqref{a88sgl} gets truncated from below and from above so that the
multiple sum is finite. By a polynomial argument, we can replace $q^M$
by $bc/aq$. If we then perform the substitution
$d\mapsto dq^{-|\mathbf k|}$ and replace $k_i$ by $m_i$, $i=1,\dots,r$,
we obtain an $A_r$ extension of Jackson's terminating balanced
very-well-poised ${}_8\phi_7$ summation (cf.\ \cite[Thm.~6.14]{M1}) which,
via a polynomial argument, is equivalent to Proposition~\ref{ar87}.

\item
If, in \eqref{a88sgl}, we let $M\to\infty$ and perform the substitution
$d\mapsto dq^{-|\mathbf k|}$, we can repeatedly apply analytic
continuation to replace $q^{k_i}$ by $a/c_ie_i$ for $i=1,\dots,r$
(in order to relax the integrality condition of the $k_i$'s),
where $e_1,\dots,e_r$ are new complex parameters. We then obtain the
$A_r$ extension of Bailey's very-well-poised ${}_6\psi_6$ summation in
Proposition~\ref{ar66}.
\end{enumerate}
\end{Remark}

\begin{Theorem}[An $A_r$ balanced very-well-poised
${}_8\psi_8$ summation formula]\label{av88s}
Let $a$, $b$, $c_1,\dots,c_r$, $d$ and $x_1,\dots,x_r$ be
indeterminate, let $k_1,\dots,k_r$ be integers, and let $M$ be a
nonnegative integer. Then
\begin{multline}\label{av88sgl}
\sum_{n_1,\dots,n_r=-\infty}^{\infty}
\frac{(1-aq^{2|\mathbf n|})}{(1-a)}\,
\frac{(b,aq^{1+M}/b,aq^{-|\mathbf k|}/C)_{|{\mathbf n}|}}
{(bq^{-M},aq/b,aq/C)_{|\mathbf n|}}\,q^{|\mathbf n|}
\prod_{1\le i<j\le r}\frac {x_iq^{n_i}-x_jq^{n_j}}{x_i-x_j}\\\times
\prod_{i,j=1}^r\frac{(c_jx_i/x_j)_{n_i}}
{(q^{1+k_j}c_jx_i/x_j)_{n_i}}
\prod_{i=1}^r\frac{(dq^{1+M}/Cx_i)_{|\mathbf n|-n_i}\,
(c_idq^{k_i}/Cx_i)_{|{\mathbf n}|}\,
(ax_iq^{-M}/d)_{n_i}}
{(dq/Cx_i)_{|\mathbf n|-n_i}\,
(c_idq^{1+M}/Cx_i)_{|{\mathbf n}|}\,
(ax_iq^{1-|\mathbf k|}/d)_{n_i}}\\
=\prod_{i,j=1}^r\frac{(qc_jx_i/c_ix_j,qx_i/x_j)_{\infty}}
{(qc_jx_i/x_j,qx_i/c_ix_j)_{\infty}}
\prod_{i=1}^r\frac{(ax_iq/c_id,c_idq/ax_i,Cx_iq/d,dq/Cx_i)_{\infty}}
{(ax_iq/d,dq/ax_i,c_idq/Cx_i,Cx_iq/c_id)_{\infty}}\\\times
\frac{(aq,q/a)_{\infty}}{(aq/C,Cq/a)_{\infty}}\,
\frac{(Cq/b,aq/bC)_M}{(aq/b,q/b)_M}\,
\prod_{i=1}^r\frac{(c_idq/Cx_i,dq/ax_i)_M}
{(c_idq/ax_i,dq/Cx_i)_M}\prod_{i,j=1}^r\frac{(qc_jx_i/x_j)_{k_i}}
{(qc_jx_i/c_ix_j)_{k_i}}\\\times
\frac{(Cq/a,q^{-M})_{|{\mathbf k}|}}
{(Cq/b,bCq^{-M}/a)_{|{\mathbf k}|}}
\prod_{i=1}^r\frac{(c_id/ax_i)_{|{\mathbf k}|}\,
(bc_id/aCx_i,c_idq^{1+M}/bCx_i)_{k_i}}
{(d/ax_i)_{|{\mathbf k}|}\,(c_id/Cx_i,c_idq^{1+M}/ax_i)_{k_i}}.
\end{multline}
\end{Theorem}

\begin{proof}
We combine the multilateral matrix inverse in Theorem~\ref{arvmi}
with the $A_r$ extension of Jackson's terminating balanced
very-well-poised ${}_8\phi_7$ summation in Proposition~\ref{ar87},
using the inverse relations \eqref{invrel}.

In particular, we have \eqref{invrel2}, with
$(g_{\mathbf k\mathbf l})_{\mathbf k,\mathbf l\in\mathbb Z^r}$
as in \eqref{garvmi},
\begin{multline*}
a_{\mathbf k}=\frac{(bq/d,bq/ad)_M}{(bCq/ad,bq/Cd)_M}
\prod_{i=1}^r\frac{(bc_iq/ax_i,bq/Cx_i)_M}{(bc_iq/Cx_i,bq/ax_i)_M}
\prod_{i,j=1}^r(c_jx_i/x_j)_{k_i}\\\times
\prod_{1\le i<j\le r}\frac{x_iq^{k_i}-x_jq^{k_j}}{x_i-x_j}
\prod_{i=1}^r\frac{(bq^{1+M}/Cx_i)_{|\mathbf k|-k_i}\,(ax_iq^{-M}/b)_{k_i}}
{(b/Cx_i)_{|\mathbf k|-k_i}\,(bc_iq^{1+M}/Cx_i)_{|\mathbf k|}}\\\times
\frac{(1-aq^{2|\mathbf k|})}{(1-a)}\,
\frac{(ad/b,bq^{1+M}/d)_{|\mathbf k|}}{(adq^{-M}/b,bq/d)_{|\mathbf k|}}\,
q^{|\mathbf k|},
\end{multline*}
and 
\begin{multline*}
b_{\mathbf l}=\prod_{1\le i<j\le r}
\frac{c_iq^{l_i}/x_i-c_jq^{l_j}/x_j}{c_i/x_i-c_j/x_j}
\prod_{i=1}^r\frac{1-bc_iq^{l_i+|\mathbf l|}/ax_i}{1-bc_i/ax_i}
\prod_{i,j=1}^r\frac 1{(qc_ix_j/c_jx_i)_{l_i}}\\\times
\frac{(q^{-M})_{|\mathbf l|}}{(bCq/ad,Cdq^{-M}/b)_{|\mathbf l|}}
\prod_{i=1}^r\frac{(bc_i/ax_i)_{|\mathbf l|}\,
(c_id/Cx_i,b^2c_iq^{1+M}/aCdx_i)_{l_i}}
{(bc_iq^{1+M}/ax_i)_{l_i}}\,x_i^{|\mathbf l|}\\\times
(-1)^{(r-1)|\mathbf l|}a^{(r-1)|\mathbf l|}b^{-r|\mathbf l|}
C^{|\mathbf l|}\,q^{|\mathbf l|+(1-r)\binom{|\mathbf l|}2},
\end{multline*}
by the $k_i\mapsto l_i$, $a\mapsto b/a$,
$c_i\mapsto c_iq^{k_i}$, $d\mapsto q^{1-|\mathbf k|}/a$,
$x_i\mapsto c_i/x_i$, $i=1,\dots,r$, case of Proposition~\ref{ar87}.
Therefore we must have \eqref{invrel1}, with
$(f_{\mathbf n\mathbf k})_{\mathbf n,\mathbf k\in\mathbb Z^r}$
as in \eqref{farvmi}, and the above sequences $a_{\mathbf k}$
and $b_{\mathbf n}$.
After simplifications and the simultaneous substitutions
$b\mapsto d$, $d\mapsto bd/a$, $k_i\mapsto n_i$, $n_i\mapsto k_i$,
$i=1,\dots,r$,
we arrive at \eqref{av88sgl}.
\end{proof}

\begin{Remark}\rm
Two special cases of Theorem~\ref{av88s} are of particular interest:
\begin{enumerate}
\item
If $c_i=q^{-k_i}$, for $i=1,\dots,r$, then the multilateral series in
\eqref{a88sgl} gets truncated from below and from above so that the
multiple sum is finite. By a polynomial argument, we can replace $q^M$
by $bc/aq$. If we then perform the substitution
$d\mapsto dq^{-|\mathbf k|}$ and replace $k_i$ by $m_i$, $i=1,\dots,r$,
we obtain an $A_r$ extension of Jackson's terminating balanced
very-well-poised ${}_8\phi_7$ summation (cf.\ \cite[Thm.~4.1]{S4}) which,
via a polynomial argument, is equivalent to Proposition~\ref{arv87}.

\item
If, in \eqref{a88sgl}, we let $b\to\infty$,
we can apply analytic continuation to replace $q^M$ by $bd/a$.
If we then perform the substitutions $d\mapsto dq^{-|\mathbf k|}$
and $c_i\mapsto c_iq^{-k_i}$, for $i=1,\dots,r$,
we can repeatedly apply analytic continuation to replace $q^{k_i}$
by $c_i/e_i$, for $i=1,\dots,r$. After subsequent relabelling
of parameters $b\mapsto d$, $c_i\mapsto a/e_i$, $d\mapsto b$,
$e_i\mapsto c_i$, for $i=1,\dots,r$, we obtain exactly the $A_r$
extension of Bailey's
very-well-poised ${}_6\psi_6$ summation in Proposition~\ref{arv66}.
\end{enumerate}
\end{Remark}

\begin{Theorem}[A $C_r$ balanced very-well-poised
${}_8\psi_8$ summation formula]\label{c88s}
Let $a$, $b$, $c_1,\dots,c_r$, $d$ and $x_1,\dots,x_r$ be
indeterminate, let $k_1,\dots,k_r$ be integers, and let $M$
be a nonnegative integer. Then
\begin{multline}\label{c88sgl}
\sum_{n_1,\dots,n_r=-\infty}^{\infty}
q^{|{\mathbf n}|}
\prod_{1\le i<j\le r}\frac {x_iq^{n_i}-x_jq^{n_j}}{x_i-x_j}
\prod_{1\le i\le j\le r}\frac{1-ax_ix_jq^{n_i+n_j}}{1-ax_ix_j}\\\times
\prod_{i,j=1}^r\frac{(c_jx_i/x_j,ax_ix_jq^{-k_j}/c_j)_{n_i}}
{(ax_ix_jq/c_i,q^{1+k_j}c_jx_i/x_j)_{n_i}}
\prod_{i=1}^r\frac{(bx_i,dx_iq^{|\mathbf k|},
ax_iq^{1+M}/b,ax_iq^{-M}/d)_{n_i}}
{(ax_iq/b,ax_iq^{1-|\mathbf k|}/d,bx_iq^{-M},dx_iq^{1+M})_{n_i}}\\
=\frac{\prod_{i,j=1}^r(qc_jx_i/c_ix_j,qx_i/x_j)_{\infty}
\prod_{1\le i\le j\le r}(ax_ix_jq,q/ax_ix_j)_{\infty}}
{\prod_{i,j=1}^r(qc_jx_i/x_j,qx_i/c_ix_j,
ax_ix_jq/c_i,c_jq/ax_ix_j)_{\infty}}\\\times
\prod_{1\le i<j\le r}
(ax_ix_jq/c_ic_j,c_ic_jq/ax_ix_j)_{\infty}
\prod_{i=1}^r\frac{(ax_iq/c_id,c_idq/ax_i,dx_iq/c_i,c_iq/dx_i)_{\infty}}
{(ax_iq/d,dq/ax_i,dx_iq,q/dx_i)_{\infty}}\\\times
\prod_{i=1}^r\frac{(ax_iq/bc_i,c_iq/bx_i,dx_iq,dq/ax_i)_M\,
(c_id/ax_i)_{|{\mathbf k}|}}
{(c_idq/ax_i,dx_iq/c_i,q/bx_i,ax_iq/b)_M\,
(dx_i,d/ax_i)_{|{\mathbf k}|}}
\prod_{1\le i<j\le r}(c_ic_j/ax_ix_j)_{k_i+k_j}^{-1}\\\times
\frac{(bd/a,dq^{1+M}/b,q^{-M})_{|\mathbf k|}\,
\prod_{i=1}^r(dx_i/c_i)_{|\mathbf k|-k_i}}
{\prod_{i=1}^r(c_iq/bx_i,bc_iq^{-M}/ax_i,c_idq^{1+M}/ax_i)_{k_i}}
\prod_{i,j=1}^r\frac{(qc_jx_i/x_j,c_jq/ax_ix_j)_{k_j}}
{(qc_jx_i/c_ix_j)_{k_j}}.
\end{multline}
\end{Theorem}

\begin{proof}
We combine the multilateral matrix inverse in Theorem~\ref{crmi}
with the $D_r$ extension of Jackson's terminating balanced
very-well-poised ${}_8\phi_7$ summation in Proposition~\ref{dr87},
using the inverse relations \eqref{invrel}.

In particular, we have \eqref{invrel2}, with
$(g_{\mathbf k\mathbf l})_{\mathbf k,\mathbf l\in\mathbb Z^r}$
as in \eqref{gcrmi},
\begin{multline*}
a_{\mathbf k}=
\prod_{i=1}^r\frac{(bc_iq/ax_i,bx_iq/c_i,
bq/adx_i,bx_iq/d)_M\,(adx_i/b,bx_iq^{1+M}/d,ax_iq^{-M}/b)_{k_i}}
{(bx_iq/c_id,bc_iq/adx_i,bx_iq,bq/ax_i)_M\,
(adx_iq^{1-M}/b,bx_iq^{1+M},bx_iq/d)_{k_i}}\\\times
q^{|\mathbf k|}\,\prod_{1\le i<j\le r}\frac{x_iq^{k_i}-x_jq^{k_j}}{x_i-x_j}
\prod_{1\le i\le j\le r}\frac{1-ax_ix_jq^{k_i+k_j}}{1-ax_ix_j}
\prod_{i,j=1}^r\frac {(c_jx_i/x_j)_{k_i}}{(ax_ix_jq/c_j)_{k_i}},
\end{multline*}
and 
\begin{multline*}
b_{\mathbf l}=\prod_{1\le i<j\le r}
\frac{c_iq^{l_i}/x_i-c_jq^{l_j}/x_j}{c_i/x_i-c_j/x_j}
\prod_{i=1}^r\frac{1-bc_iq^{l_i+|\mathbf l|}/ax_i}{1-bc_i/ax_i}
\prod_{i,j=1}^r\frac 1{(qc_ix_j/c_jx_i)_{l_i}}\\\times
\prod_{1\le i<j\le r}\frac 1{(c_ic_j/ax_ix_j)_{l_i+l_j}}
\prod_{i=1}^r(bc_i/ax_i)_{|\mathbf l|}(bx_iq/c_i)_{|\mathbf l|-l_i}
c_i^{rl_i}x_i^{-rl_i}\\\times
\frac{(d,b^2q^{1+M}/ad,q^{-M})_{|\mathbf l|}}
{\prod_{i=1}^r(bc_iq/adx_i,dc_iq^{-M}/bx_i,bc_iq^{1+M}/ax_i)_{l_i}}\,
b^{-r|\mathbf l|}q^{-r\binom{|\mathbf l|}2+r\sum_{i=1}^r\binom{l_i+1}2},
\end{multline*}
by the $k_i\mapsto l_i$, $a\mapsto b/a$, $b\mapsto d$,
$c_i\mapsto c_iq^{k_i}$, $cd\mapsto 1/a$, $x_i\mapsto c_i/x_i$,
$i=1,\dots,r$, case of Proposition~\ref{dr87}.
Therefore we must have \eqref{invrel1}, with
$(f_{\mathbf n\mathbf k})_{\mathbf n,\mathbf k\in\mathbb Z^r}$
as in \eqref{fcrmi}, and the above sequences $a_{\mathbf k}$
and $b_{\mathbf n}$.
After simplifications and the simultaneous substitutions
$b\mapsto d$, $d\mapsto bd/a$, $k_i\mapsto n_i$, $n_i\mapsto k_i$,
$i=1,\dots,r$,
we arrive at \eqref{c88sgl}.
\end{proof}

\begin{Remark}\rm
Two special cases of Theorem~\ref{c88s} are of particular interest:
\begin{enumerate}
\item
If $c_i=q^{-k_i}$, for $i=1,\dots,r$, then the multilateral series in
\eqref{c88sgl} gets truncated from below and from above so that the
multiple sum is finite. By a polynomial argument, we can replace $q^M$
by $bc/aq$. If we then perform the substitution
$d\mapsto dq^{-|\mathbf k|}$ and replace $k_i$ by $m_i$, $i=1,\dots,r$,
we obtain the $C_r$ extension of Jackson's terminating balanced
very-well-poised ${}_8\phi_7$ summation in Proposition~\ref{cr87}.

\item
If, in \eqref{c88sgl}, we let $M\to\infty$ and perform the substitution
$d\mapsto dq^{-|\mathbf k|}$, we can repeatedly apply analytic
continuation to replace $q^{k_i}$ by $a/c_ie_i$ for $i=1,\dots,r$
(in order to relax the integrality condition of the $k_i$'s),
where $e_1,\dots,e_r$ are new complex parameters. We then obtain the
$C_r$ extension of Bailey's very-well-poised ${}_6\psi_6$ summation in
Proposition~\ref{cr66}.
\end{enumerate}
\end{Remark}

\begin{appendix}
\section{}
\label{secdr66}

It is indeed quite interesting that concerning the multilateral matrix
inversion result in Theorem~\ref{crmi},
$(f_{\mathbf n\mathbf k})_{\mathbf n,\mathbf k\in\mathbb Z^r}$
happens to be the left-inverse of
$(g_{\mathbf k\mathbf l})_{\mathbf k,\mathbf l\in\mathbb Z^r}$,
but {\em not} the right-inverse (unless in special cases, e.g.,
when both matrices are upper- or lower-triangular).
Let us assume, for a moment, that
$(f_{\mathbf n\mathbf k})_{\mathbf n,\mathbf k\in\mathbb Z^r}$
would also be the right-inverse of 
$(g_{\mathbf k\mathbf l})_{\mathbf k,\mathbf l\in\mathbb Z^r}$
(which cannot be justified, see also Remark~\ref{rem}).
By combining these matrices with the $C_r$ very-well-poised
${}_6\psi_6$ summation formula in Proposition~\ref{cr66},
by virtue of multidimensional inverse relations
one would be able to deduce a ``new'' $D_r$ very-well-poised
${}_6\psi_6$ summation formula which, however, turns out
to be {\em false} for $r>1$, as the series does not converge.

In particular, we have \eqref{rotinv2} with
$(g_{\mathbf k\mathbf l})_{\mathbf k,\mathbf l\in\mathbb Z^r}$
as in \eqref{gcrmi},
\begin{multline*}
a_{\mathbf l}=
\prod_{1\le i<j\le r}(ax_ix_jq/c_ic_j,aq/e_ie_jx_ix_j)_{\infty}
\prod_{1\le i\le j\le r}(ax_ix_jq,q/ax_ix_j)_{\infty}\\\times
\frac{(bq/d)_{\infty}}{(a^{r+1}bq/CdE)_{\infty}}
\prod_{i,j=1}^r\frac{(ax_iq/c_ie_jx_j,qx_i/x_j)_{\infty}}
{(ax_iq/e_jx_j,q/e_jx_ix_j,ax_ix_jq/c_i,qx_i/c_ix_j)_{\infty}}\\\times
\prod_{i=1}^r\frac{(ax_iq/c_id,aq/de_ix_i,bx_iq/c_i,bq/e_ix_i)_{\infty}}
{(bx_iq,bq/ax_i,ax_iq/d,q/dx_i)_{\infty}}\\\times
\prod_{i,j=1}^r(c_ie_jx_j/ax_i)_{l_i}
\prod_{1\le i<j\le r}(c_ic_j/ax_ix_j)_{l_i+l_j}
\prod_{i=1}^r\frac{(c_id/ax_i)_{l_i}}{(bq/e_ix_i)_{|\mathbf l|}\,
(bx_iq/c_i)_{|\mathbf l|-l_i}}\\\times
\frac 1{(bq/d)_{|\mathbf l|}}\,
\left(\frac{a^rb^r}{CdE}\right)^{|\mathbf l|}
q^{(r-1)\left(\binom{|\mathbf l|}2-\sum_{i=1}^r\binom{l_i+1}2\right)}
\prod_{i=1}^r\left(\frac{x_i}{c_i}\right)^{(r-1)l_i}
\end{multline*}
and
\begin{equation*}
b_{\mathbf k}=
\prod_{i,j=1}^r\frac{(e_jx_ix_j)_{k_i}}{(ax_iq/e_jx_j)_{k_i}}
\prod_{i=1}^r\frac{(dx_i)_{k_i}}{(ax_iq/d)_{k_i}}\cdot
\left(\frac{a^rb}{CdE}\right)^{|\mathbf k|},
\end{equation*}
by the $b\mapsto aq^{-|\mathbf l|}/b$, $c_i\mapsto c_iq^{l_i}$,
$i=1,\dots,r$, case of Proposition~\ref{cr66}.
As we are (erroneously) assuming that
$(f_{\mathbf n\mathbf k})_{\mathbf n,\mathbf k\in\mathbb Z^r}$
is the right-inverse of  
$(g_{\mathbf k\mathbf l})_{\mathbf k,\mathbf l\in\mathbb Z^r}$,
we must have \eqref{rotinv1}, with
$(f_{\mathbf n\mathbf k})_{\mathbf n,\mathbf k\in\mathbb Z^r}$
as in \eqref{farmi}, and the above sequences $a_{\mathbf n}$
and $b_{\mathbf k}$. After simplifications and the simultaneous
substitutions $a\mapsto bc/a^2$, $b\mapsto bc/a$,
$c_i\mapsto aq^{-k_i}/c_i$, $d\mapsto cd/a$, 
$e_i\mapsto bcc_ie_iq^{k_i}/a^3$, $x_i\mapsto aq^{-k_i}/c_ix_i$,
$i=1,\dots,r$, we can get rid of the $k_i$ and would arrive at
the following identity:
\begin{multline}\label{dr66gl}
\sum_{n_1,\dots,n_r=-\infty}^{\infty}
\prod_{1\le i<j\le r}\frac {x_iq^{n_i}-x_jq^{n_j}}{x_i-x_j}
\prod_{i=1}^r\frac {1-ax_iq^{n_i+|\mathbf n|}}{1-ax_i}
\prod_{1\le i<j\le r}(a^2x_ix_j/bc)_{n_i+n_j}\\\times
\prod_{i,j=1}^r\frac{(e_jx_i/x_j)_{n_i}}
{(ax_iq/c_jx_j,ac_ix_ix_jq/bc)_{n_i}}
\prod_{i=1}^r\frac{(bc/c_ix_i,c_ix_i)_{|\mathbf n|}\,(dx_i)_{n_i}}
{(ax_iq/e_i)_{|\mathbf n|}\,(bc/ax_i)_{|\mathbf n|-n_i}}\,x_i^{n_i}\\\times
\frac 1{(aq/d)_{|\mathbf n|}}
\left(\frac{a^2q}{bcdE}\right)^{|{\mathbf n}|}\,
q^{-e_2(\mathbf n)}\\
=\prod_{1\le i<j\le r}\frac{(a^2x_ix_jq/bc)_{\infty}}
{(a^2x_ix_jq/bce_ie_j)_{\infty}}
\prod_{i,j=1}^r\frac{(ax_iq/c_je_ix_j,ac_ix_ix_jq/bce_j,qx_i/x_j)_{\infty}}
{(ax_iq/c_jx_j,ac_ix_ix_jq/bc,qx_i/e_ix_j)_{\infty}}\\\times
\frac{(aq/dE)_{\infty}}{(aq/d)_{\infty}}
\prod_{i=1}^r\frac{(ax_iq,q/ax_i,ax_iq/bc,aq/c_idx_i,ac_ix_iq/bcd)_{\infty}}
{(ax_iq/e_i,c_ix_iq/bc,q/c_ix_i,q/dx_i,a^2x_iq/bcde_i)_{\infty}},
\end{multline}
where $e_2(\mathbf n)$ is the second elementary symmetric function
of $(n_1,\dots,n_r)$.

Now, due to the factor $q^{-e_2(\mathbf n)}$ appearing in the summand,
the series in \eqref{dr66gl} does {\em not} converge for $r\ge 2$.
Therefore, the identity as stated is false. However, it is valid
for $r\ge 2$ whenever the series terminates. For instance, when
$c_i=a$ and $e_i=q^{-m_i}$, for $i=1,\dots,r$, \eqref{dr66gl} reduces
to Bhatnagar's $D_r$ terminating very-well-poised ${}_6\phi_5$
summation, derived in \cite[Thm.~2]{Bh}.
\end{appendix}

\end{document}